\documentclass[12pt,reqno]{amsart}
\usepackage{amssymb, amsmath}
\usepackage{amsfonts}
\usepackage[margin=1.0in]{geometry}
\usepackage{esint}
\usepackage{enumitem}
\usepackage{latexsym}

\allowdisplaybreaks 
\theoremstyle{plain}

\newcommand{\C}{\mathbb C}

\newtheorem{thm}{Theorem}[section]
\newtheorem{prop}[thm]{Proposition}
\newtheorem{lem}[thm]{Lemma}
\newtheorem{cor}[thm]{Corollary}
\newtheorem{defn}[thm]{Definition}
\newtheorem{fact}[thm]{Fact}

\newtheorem{example}[thm]{Example}

\newtheorem{rem}[thm]{Remark}

\makeatletter\@addtoreset{equation}{section} \makeatother

\begin{document}

\title{Complex solutions to Maxwell's equations}
\author{Sachin Munshi$^*$ and Rongwei Yang}
\thanks{
$^*$ Corresponding author}
\address[Sachin Munshi]{Department of Mathematics and Statistics, SUNY at Albany, Albany, NY 12222, U.S.A.}
\email[Sachin Munshi]{sacmun86@gmail.com, smunshi@albany.edu}

\address[Rongwei Yang]{Department of Mathematics and Statistics, SUNY at Albany, Albany, NY 12222, U.S.A. }
\email[Rongwei Yang]{ryang@albany.edu}


\date{}

\maketitle

\begin{abstract}
This paper provides a view of Maxwell's equations from the perspective of complex variables. The study is made through complex differential forms and the Hodge star operator in $\C^2$ with respect to the Euclidean and the Minkowski metrics. It shows that holomorphic functions give rise to nontrivial solutions, and the inner product between the electric and the magnetic fields is considered in this case.
Further, it obtains a simple necessary and sufficient condition regarding harmonic solutions to the equations. In the end, the paper gives an interpretation of the Lorenz gauge condition in terms of the codifferential operator. 
\end{abstract}

\section{Introduction}
Named after the physicist and mathematician James C. Maxwell, Maxwell's
equations form the foundation of classical electromagnetism, optics,
and electrodynamics (see for instance \cite{baez1994gauge, garrity2015electricity, maxwell1865viii}). 
They are a set of partial differential equations that
describe the interactions between the electric and magnetic fields that emerge
from distributions of electric charges and currents, and of course, how these
fields change in time. Maxwell's equations are the differential equations
\begin{subequations}

\begin{eqnarray}
\nabla\cdotp\mathbf{B} & = & 0,\label{eq:1.1a}\\
\nabla\times\mathbf{E}+\frac{\partial\mathbf{B}}{\partial t} & = & 0,\label{eq:1.1b}\\
\nabla\cdotp\mathbf{E} & = & \frac{\rho}{\epsilon_{0}},\label{eq:1.1c}\\
\nabla\times\mathbf{B}-\frac{1}{c^{2}} \frac{\partial\mathbf{E}}{\partial t} & = & \mu_{0} \mathbf{J},\label{eq:1.1d}
\end{eqnarray}

\end{subequations} 

\noindent where $\mathbf{E}=(E_1, E_2, E_3)$ is the \textit{electric field}, $\mathbf{B}=(B_1, B_2, B_3)$ is the \textit{magnetic field}, the scalar $\rho$ is the \textit{electric charge density}, and the vector $\mathbf{J}$ is the \textit{electric current density vector}. Moreover, $\epsilon_{0}$ is the \textit{vacuum permittivity}, $\mu_{0}$ is the \textit{vacuum permeability}, and $c:= 1/ \sqrt{\epsilon_{0}\mu_{0}}$ is the \textit{speed of light} in vacuum. For ease of use, most modern physicists and mathematicians simply set $c=1$, as shall we throughout this paper. For convenience, we fix the notation $x_0=ct=t$. Then in a more compact form, equations (\ref{eq:1.1a})-(\ref{eq:1.1d}) can be written in differential forms over Minkowski space-time as
\begin{eqnarray}
dF & = & 0,\label{eq:1.2}\\
\star \text{ } d \star F & = & J,\label{eq:1.3}
\end{eqnarray}
where \begin{align}
F & =-dx_{0}\wedge\left(E_{1}dx_{1}+E_{2}dx_{2}+E_{3}dx_{3}\right)\label{eq:3.17}\\
 & -B_{1}dx_{2}\wedge dx_{3}+B_{2}dx_{1}\wedge dx_{3}-B_{3}dx_{1}\wedge dx_{2},\nonumber 
\end{align}
and it is referred to as the \textit{Faraday} $2$-\textit{form}. Here, $J$ is the \textit{current} $1$-\textit{form}, $d$ is the \textit{exterior differential operator}, and $\star$ is the \textit{Hodge star operator} (\cite{baez1994gauge}). These will be described in detail in this paper. By the Poincar\'{e} lemma, Equation (\ref{eq:1.2}) implies that locally $F=d\omega$ for some differentiable $1$-form
$\omega=\eta_0dx_0+\eta_1dx_1+\eta_2dx_2+\eta_3dx_3.$ This paper investigates the case when $\eta_j, 0\leq j\leq 3$, are all harmonic functions in $(x_0, x_1, x_2, x_3)$ with respect to the Euclidean or Minkowski metric. If one takes advantage of the identification $\mathbb{R}^{4}\simeq\mathbb{C}^{2}$,
i.e. 
\[
\mathbb{R}^{4}\ni\left(x_{0},x_{1},x_{2},x_{3}\right)\mapsto\left(x_{0}+ix_{1},x_{2}+ix_{3}\right):=\left(z_{1},z_{2}\right)\in\mathbb{C}^{2},
\]
then $\omega$ can be written as $\omega(z)=f_1dz_1+f_2dz_2+f_{\bar{1}}d\bar{z}_{1}+f_{\bar{2}}d\bar{z}_{2}$. The following is the main result.
\begin{thm}
Let $f_{j},f_{\bar{j}},\text{ }j=1,2$ be harmonic functions on $\mathbb{C}^{2}$. Then the complex differential form $F_\omega:=d\omega$ is a solution to the source-free Maxwell's equations in the Euclidean metric if and only if $\bar{\partial}_{1}f_{1}+\bar{\partial}_{2}f_{2}+\partial_{1}f_{\bar{1}}+\partial_{2}f_{\bar {2}}$ is constant.
\end{thm}
A parallel theorem holds with respect to the Minkowski metric on $\mathbb{R}^{4}\simeq\mathbb{C}^{2}$. A solution $F_\omega$ to the Maxwell's equations with respect to the Minkowski metric is said to be {\em wavelike} if each of the functions $f_{j},f_{\bar{j}},\text{ }j=1,2$, is a solution to the d'Alembertian equation (or wave equation)
\[\left(\frac{\partial^{2}}{\partial x_{0}^{2}}- \frac{\partial^{2}}{\partial x_{1}^{2}}- \frac{\partial^{2}}{\partial x_{2}^{2}}- \frac{\partial^{2}}{\partial x_{3}^{2}}\right)u=0.\]
Known solutions to the Maxwell's equations are all wavelike, so the next result seems surprising.
\begin{thm}
There exist non-wavelike solutions to the source-free Maxwell's equations with respect to the Minkowski metric.
\end{thm}

The paper is organized as follows.

\tableofcontents

\section{Preliminaries}

Recall that $\mathbb{R}^{4}$ is just Euclidean space of real dimension
4. As a vector space, a point in $\mathbb{R}^{4}$ may be considered
as a row vector $\mathbf{x}=\left(x_{0},x_{1},x_{2},x_{3}\right),x_{i}\in\mathbb{R}$.
Note that in the Lorentzian signature, $\mathbb{R}^{4}$ is denoted
as $\mathbb{R}^{1,3}$, and is referred to as Minkowski space-time,
with $x_{0}$ as the time variable (denoted $ct$), where $c$ is the speed of light, and $x_{1},x_{2},x_{3}$ are the space variables. Going back to viewing $\mathbb{R}^{4}$ in
the Euclidean metric, we have the identification $\mathbb{R}^{4}\simeq\mathbb{C}^{2}$,
i.e. 
\[
\mathbb{R}^{4}\ni\left(x_{0},x_{1},x_{2},x_{3}\right)\mapsto\left(x_{0}+ix_{1},x_{2}+ix_{3}\right):=\left(z_{1},z_{2}\right)\in\mathbb{C}^{2}.
\]

\subsection{Complex Differential Forms}

We first introduce some preliminaries on complex differential forms. Consider $\mathbb{C}^{n}$ with points given by coordinates $z=\left(z_{1},z_{2},\dots,z_{n}\right)$.
The \textit{tangent space} of $\mathbb{C}^{n}$ is $$T\left(\mathbb{C}^{n}\right)=\text{span}\left\{ \frac{\partial}{\partial z_{k}},\frac{\partial}{\partial\bar{z}_{k}}:1\leq k\leq n\right\}, $$ and 
the \textit{cotangent space} is given by $$T^{\ast}\left(\mathbb{C}^{n}\right)=\text{span}\left\{ dz_{k},d\bar{z}_{k}:1\leq k\leq n\right\}. $$ 
\begin{defn}
Let $f$ be a smooth function on a domain $\mathcal{M}\subset\mathbb{C}^{n}$.
Consider the linear operators $\partial,\bar{\partial},d$ defined
to act on $f$ as follows:
\begin{align*}
\partial f  =\sum_{k=1}^{n}\frac{\partial f}{\partial z_{k}}dz_{k},\hspace{0.5cm}
\bar{\partial}f =\sum_{k=1}^{n}\frac{\partial f}{\partial\bar{z}_{k}}d\bar{z}_{k},\hspace{0.5cm}
df =\left(\partial+\bar{\partial}\right)f.
\end{align*}
Here, $d$ is called the \textit{complex exterior differential operator},
while $\partial,\bar{\partial}$ are called the \textit{Dolbeault operators}. 
\end{defn}

Recall that a smooth function $f$ on a domain $\mathcal{M}$
is said to be {\it holomorphic} if $\bar{\partial} f=0$ everywhere on $\mathcal{M}$.
Clearly, this means that $f$ is analytic in each variable (\cite{range2013holomorphic}). 
Note that if $f$ is holomorphic, then $df=\partial f$. The following fact is well-known.
\begin{fact}
$\partial^{2}=\bar{\partial}^{2}=d^{2}=0$.
\end{fact}
For a multi-index $I=i_1i_2\cdots i_p$ we assume $i_1<i_2<\cdots <i_p$ and define its length $|I|=p$.
\begin{defn}
The space of complex $\left(p,q\right)$-forms on $\mathcal{M}$ is defined as
$$\Omega^{p,q}\left(\mathcal{M}\right):=\left\{ \sum_{|I|=p,|J|=q}f_{I,J}dz_{I}\wedge d\bar{z}_{J}: p+q\leq 2n, f_{I,J}\in C^{\infty}\left(\mathcal{M}\right)\right\}.$$
\end{defn}
We may drop the $\mathcal{M}$ from the definition for convenience. Clearly, $\Omega^{1,0}$ is the space of complex differential forms containing only the $dz_{k}$ terms, and $\Omega^{0,1}$
is the space of forms containing only the $d\bar{z}_{k}$ terms.
Then in terms of the exterior product on differential forms we have 
\begin{equation*}
\Omega^{p,q}= \underbrace{\Omega^{1,0} \wedge \cdots \wedge \Omega^{1,0}}_{p} \wedge
\underbrace{\Omega^{0,1} \wedge \cdots \wedge \Omega^{0,1}}_{q}.
\end{equation*}
Slightly abusing notation, we shall adhere to the following definition throughout this paper.
\begin{defn}
$\Omega^{k}:=\bigoplus_{p+q=k}\Omega^{p,q}$ is the space of
all \textit{complex differential forms of total degree} $k=p+q$.
\end{defn}

\begin{defn}
For each $p>0$, the forms given by $\sum_{|I|=p}f_{I}dz_{I}$, where each
 $f_{I}$ is holomorphic, are called \textit{holomorphic p-forms},
and they form a holomorphic section of $\Omega^{p,0}$.
\end{defn}

Note that if $\eta=\sum_{|I|=p}f_{I}d\bar{z}_{I}$, with $f_{I}$
holomorphic, then $\bar{\partial}\eta=0$. 

\subsection{Hodge Star Operator $\star$}

We briefly go over some basics of Hodge theory withholding any discussion
on topology or manifold theory. We identify $\C^n$ with ${\mathbb R}^{2n}$ with the representation
$z_k=x_{2k-2}+ix_{2k-1},\ k=1, 2, ..., n$. Given a nondegenerate self-adjoint $2n\times 2n$ matrix $g=\left(g_{st}\right)$, it induces a  sesquilinear form $\langle \cdot, \cdot \rangle$ on the tangent space $T\left(\mathbb{R}^{2n}\right)$ with complex coefficients by the evaluations
\[\left\langle \alpha \frac{\partial}{\partial x_{s}}, \beta \frac{\partial}{\partial x_{t}}\right\rangle_g:=\alpha \overline{\beta}g_{st},\]
where $\alpha$ and $\beta$ are complex numbers and $0\leq s, t\leq 2n-1$. Then on the cotangent space $T^*\left(\mathbb{R}^{2n}\right)$ one has the corresponding sesquilinear form given by
\[\langle \alpha dx_s, \beta dx_t\rangle_g:=\alpha \overline{\beta}g^{st},\]
where $(g^{st})=g^{-1}$. Using the fact $z_k=x_{2k-2}+ix_{2k-1}$ above and the following representation
\[\frac{\partial}{\partial z_{k}}=\frac{1}{2}\left(\frac{\partial}{\partial x_{2k-2}}-i\frac{\partial}{\partial x_{2k-1}}\right), k= 1, 2, ..., n,\]
one may regard the aforementioned sesquilinear form $\langle \cdot, \cdot \rangle_g$ as a sesquilinear form on $T^*\left(\mathbb{C}^{n}\right)$. Further, it can be extended to a sesquilinear form on $\Omega^{p}$ such that for $\eta=\eta_1\wedge\cdots \wedge \eta_p,\ \xi=\xi_1\wedge\cdots \wedge \xi_p$ one has  \[
\langle \eta,\xi \rangle_g= \det [\langle \eta_{s},\xi_{t}\rangle_g]_{s,t=1}^{p}, \text{ }\eta, \xi \in \Omega^{p}.
\]
One observes that if the matrix $g$ is positive definite then $\langle \cdot, \cdot \rangle_g$ is an inner product on the set of constant $p$-forms for each $1\leq p\leq 2n$.

\begin{defn}\label{de:hodge}
The \textit{Hodge star operator} $\star:\Omega^{p}\left(\mathbb{C}^{n}\right)\rightarrow\Omega^{2n-p}\left(\mathbb{C}^{n}\right)$ with respect to the bilinear form $\langle \cdot, \cdot \rangle_g$ 
is a linear operator such that for $\eta, \xi \in \Omega^{p}\left(\C^n\right)$ one has
\[\eta\wedge \star \bar{\xi}=\langle \eta, \xi \rangle_g \textup{vol}_g,\]
where $\textup{vol}_g=\left(\frac{i}{2}\right)^n\sqrt{|\det g|}dz_1\wedge d\bar{z}_{1}\wedge \cdots \wedge dz_n\wedge d\bar{z}_{n}$ is the volume form \textup{(\cite{huybrechts2006complex})}.
\end{defn}

\begin{example}\label{ex:hodge}
For the Euclidean metric on ${\mathbb R}^{2n}$ we have $g=I_{2n}$. Let $\Omega^{p}\left(\mathbb{R}^{2n}\right)$
denote the space of real differential $p$-forms over $\mathbb{R}^{2n}$. Then it is well-known that
\begin{equation}
\star\left(dx_{0}\wedge dx_{2}\wedge\cdots\wedge dx_{p-1}\right)=dx_{p}\wedge dx_{p+2}\wedge\cdots\wedge dx_{2n-1}.\label{eq:2.1}
\end{equation}

And for a permutation $\sigma=\left(i_{0},i_{1},\dots,i_{2n-1}\right)$,
we have 
\[
\star\left(dx_{i_{0}}\wedge dx_{i_{1}}\wedge\cdots\wedge dx_{i_{p-1}}\right)=\left(-1\right)^{\sigma}dx_{i_{p}}\wedge dx_{i_{p+1}}\wedge\cdots\wedge dx_{i_{2n-1}}
\]
depending on the parity of $\sigma$.  In particular, $\star\left(dx_{0}\wedge dx_{1}\wedge\cdots\wedge dx_{2n-1}\right)=1$
\footnote{$dx_{0}\wedge dx_{1}\wedge\cdots\wedge dx_{2n-1}$ is the volume
form.}. Moreover, for a $p$-form $\omega$,
\[\star^{2}\omega=\left(-1\right)^{p\left(2n-p\right)}\omega.\] This is only true in the Euclidean metric. Generally, $\star^{2}\omega=\left(-1\right)^{p\left(2n-p\right)}s\omega$,
where $s$ is the parity of the signature of the inner product defined by the metric. 
\end{example}

The next two examples give the Hodge dual on $\C^2$ with respect to the Euclidean metric and the Minkowski metric, respectively, and they will be used later. Using the earlier
identification $z_{1}=x_{0}+ix_{1},z_{2}=x_{2}+ix_{3}$, we have 
\[
dz_{1}=dx_{0}+idx_{1},dz_{2}=dx_{2}+idx_{3},d\bar{z}_{1}=dx_{0}-idx_{1},d\bar{z}_{2}=dx_{2}-idx_{3},
\]
and these complex $\left(1,0\right)$- and $\left(0,1\right)$-forms
span the complex cotangent space $T^{\ast}\left(\mathbb{C}^{2}\right)$. Clearly, $T^{\ast}\left(\mathbb{C}^{2}\right)$
is also spanned by the real forms $dx_{k}, 0\leq k\leq 3$ with complex coefficients.

\begin{example}\label{duals}
One may directly use \textup{Definition \ref{de:hodge}} to compute the Hodge dual of complex differential forms, or one may first write them as real forms, apply \textup{Example \ref{ex:hodge}} and then convert
back to complex forms. We shall use the latter approach.

Let's start with $\star dz_{1}$. Since $dz_{1}=dx_{0}+idx_{1}$, it follows that
\begin{align*}
\star dz_{1} & =\star dx_{0}+i\star dx_{1}\\
 & =\left(dx_{1}\wedge dx_{2}\wedge dx_{3}\right)-i\left(dx_{0}\wedge dx_{2}\wedge dx_{3}\right)\\
 & =\left(dx_{1}-idx_{0}\right)\wedge\left(dx_{2}\wedge dx_{3}\right)\\
 & =\left(dx_{1}-idx_{0}\right)\wedge\frac{i}{2}\left(dz_{2}\wedge d\bar{z}_{2}\right)\\
 & =\left(dx_{0}+idx_{1}\right)\wedge\frac{1}{2}\left(dz_{2}\wedge d\bar{z}_{2}\right)\\
 & =\frac{1}{2}\left(dz_{1}\wedge dz_{2}\wedge d\bar{z}_{2}\right),
\end{align*}
where the fourth line follows from the fact that
\begin{align*}
dz_{2}\wedge d\bar{z}_{2} & =\left(dx_{2}+idx_{3}\right)\wedge\left(dx_{2}-idx_{3}\right)=-2idx_{2}\wedge dx_{3}.
\end{align*}

Likewise one verifies that
\begin{align*}
\star dz_{2} =-\frac{1}{2}\left(dz_{1}\wedge dz_{2}\wedge d\bar{z}_{1}\right),
\star d\bar{z}_{1} =\frac{1}{2}\left(dz_{2}\wedge d\bar{z}_{1}\wedge d\bar{z}_{2}\right),
\star d\bar{z}_{2} =-\frac{1}{2}\left(dz_{1}\wedge d\bar{z}_{1}\wedge d\bar{z}_{2}\right),
\end{align*}
and 
\begin{align*}
&\star\left(dz_{1}\wedge d\bar{z}_{1}\right) =dz_{2}\wedge d\bar{z}_{2}, 
&\star\left(dz_{2}\wedge d\bar{z}_{2}\right) =dz_{1}\wedge d\bar{z}_{1},\hspace{5mm}
&\star\left(dz_{1}\wedge d\bar{z}_{2}\right) =-dz_{1}\wedge d\bar{z}_{2},\\
&\star\left(dz_{2}\wedge d\bar{z}_{1}\right) =-dz_{2}\wedge d\bar{z}_{1},
&\star\left(dz_{1}\wedge dz_{2}\right) =dz_{1}\wedge dz_{2},  \hspace{5mm}
&\star\left(d\bar{z}_{1}\wedge d\bar{z}_{2}\right) =d\bar{z}_{1}\wedge d\bar{z}_{2}.
\end{align*}
The Hodge star of $3$-forms can be computed using the fact that $\star^2 \omega=(-1)^{p(n-p)} \omega$, and one has
\begin{align*}
&\star\left(dz_{1}\wedge dz_{2}\wedge d\bar{z}_{1}\right) =2dz_{2},\hspace{5mm}
\star\left(dz_{1}\wedge dz_{2}\wedge d\bar{z}_{2}\right) =-2dz_{1},\\
&\star\left(dz_{1}\wedge d\bar{z}_{1}\wedge d\bar{z}_{2}\right) =2d\bar{z}_{2},\hspace{5mm}
\star\left(dz_{2}\wedge d\bar{z}_{1}\wedge d\bar{z}_{2}\right) =-2d\bar{z}_{1}.
\end{align*}
 Finally, for the unique $\left(2,2\right)$-form, we have 
\begin{align*}
\star\left(dz_{1}\wedge dz_{2}\wedge d\bar{z}_{1}\wedge d\bar{z}_{2}\right)=4.\\
\end{align*}
\end{example}

\begin{example}\label{mink}
The Minkowski metric on ${\mathbb R}^{1,3}$ has the signature $\left(+ - - -\right)$ represented by the metric matrix
\begin{equation}
g_{\textbf{mink}}=\begin{pmatrix}
1 & 0 & 0 & 0\\
0 & -1 & 0 & 0\\
0 & 0 & -1 & 0\\
0 & 0 & 0 & -1 
\end{pmatrix}.
\label{eq:3.19a}
\end{equation}
In this case, the inner product of complex $1$-forms is determined by the facts:
\[\langle dz_1, dz_1\rangle =\langle dz_1, dz_2\rangle=\langle dz_1, d\bar{z}_2\rangle=\langle dz_2, d\bar{z}_2\rangle=0\]
and \[\langle dz_1, d\bar{z}_1\rangle=\langle dz_2, dz_2\rangle=2.\]
Calculation of the Hodge star operator on $\mathbb{C}^2$ with respect to the Minkowski metric is similar to that in \textup{Example \ref{duals}}. Here we only list its action on 2-forms and 3-forms for later use:
\begin{align*}
&\star\left(dz_{1}\wedge dz_{2}\right) = -dz_{2}\wedge d\bar{z}_{1}, & \star\left(dz_{1}\wedge d\bar{z}_{2}\right) = -d\bar{z}_{1}\wedge d\bar{z}_{2}, & & \star\left(dz_{2}\wedge d\bar{z}_{1}\right) = dz_{1}\wedge dz_{2},\\ &\star\left(d\bar{z}_{1}\wedge d\bar{z}_{2}\right) = dz_{1}\wedge d\bar{z}_{2}, & \star\left(dz_{1}\wedge d\bar{z}_{1}\right) = -dz_{2}\wedge d\bar{z}_{2}, &  & \star\left(dz_{2}\wedge d\bar{z}_{2}\right) = dz_{1}\wedge d\bar{z}_{1};
\end{align*}
and 
\begin{align*}
&\star\left(dz_{1}\wedge dz_{2}\wedge d\bar{z}_1\right) =2dz_2, \hspace{1cm}  \star\left(dz_{1}\wedge dz_{2}\wedge d\bar{z}_2\right) =2dz_{\bar{1}},\\
&\star\left(dz_{1}\wedge d\bar{z}_1\wedge d\bar{z}_2\right) =2d\bar{z}_2, \hspace{1cm}
\star\left(dz_{2}\wedge d\bar{z}_1\wedge d\bar{z}_2\right) =2dz_{1}.
\end{align*}

\end{example}

\subsection{Self-dual and Anti-self-dual Forms}

We now focus on the Hodge star operator on the complex differential forms on $\C^2$ with respect to the Euclidean metric and the Minkowski metric.

{\bf 1}. Under the Euclidean metric, we have 
$\star^2\omega=\omega$ for every $\omega\in \Omega^{2}(\mathbb{C}^{2})$.
Self-dual and anti-self-dual forms are the eigenvectors corresponding to the eigenvalues $1$ and $-1$, respectively, of the Hodge star operator on $\C^2$.
\begin{defn}
A differential form $\omega$ is said to be \textit{self-dual} if
it is equal to its Hodge dual, i.e. $\star\omega=\omega$. If $\star\omega=-\omega$,
then $\omega$ is said to be \textit{anti-self-dual}.
\end{defn}

Out of the complex differential forms we considered in the previous subsection, six of them correspond
to the pair $\left(p,q\right)$ such that $p+q=2$. Let $\Omega_{+}^{2},\Omega_{-}^{2}$ denote
the bases of self-dual and anti-self-dual forms, respectively, in
$\Omega^{2}$. Then by the computations in Example \ref{duals} we have that 
\begin{align}
\Omega_{+}^{2} & =\text{span}\left\{ dz_{1}\wedge dz_{2},d\bar{z}_{1}\wedge d\bar{z}_{2},dz_{1}\wedge d\bar{z}_{1}+dz_{2}\wedge d\bar{z}_{2}\right\} ,\label{eq:2.5}\\
\Omega_{-}^{2} & =\text{span}\left\{ dz_{1}\wedge d\bar{z}_{2},dz_{2}\wedge d\bar{z}_{1},dz_{1}\wedge d\bar{z}_{1}-dz_{2}\wedge d\bar{z}_{2}\right\} .\nonumber 
\end{align}

{\bf 2}. Under the Minkowski metric, we have $\star^{2}\omega= -\omega$ for any $\omega\in \Omega^2(\mathbb{C}^{2})$.
Hence $\star$ has eigenvalues $i, -i$. With a bit of abuse of terminology, the eigenspaces corresponding to them are often also called self-dual and anti-self-dual forms, respectively. For consistency we shall also denote them by $\Omega^2_+$ and $\Omega^2_-$, respectively.
Then by the computations in Example \ref{mink}, we have 
\begin{align}
\Omega_{+}^{2} & =\text{span}\left\{ dz_{1}\wedge dz_{2}+idz_2\wedge d\bar{z}_{1}, d{z}_{1}\wedge d\bar{z}_{1}+idz_{2}\wedge d\bar{z}_{2}, d{z}_{1}\wedge d\bar{z}_{2}+id\bar{z}_{1}\wedge d\bar{z}_{2}\right\} ,\label{eq:2.5}\\
\Omega_{-}^{2} & =\text{span}\left\{ dz_{1}\wedge dz_{2}-idz_2\wedge d\bar{z}_{1}, d{z}_{1}\wedge d\bar{z}_{1}-idz_{2}\wedge d\bar{z}_{2}, d{z}_{1}\wedge d\bar{z}_{2}-id\bar{z}_{1}\wedge d\bar{z}_{2}\right\}\label{eq:2.5}.\nonumber 
\end{align}

\subsection{The Hodge Laplacian}

The Hodge star operator gives rise to a Hermitian bilinear form on compactly supported differentiable $p$-forms defined by
\[(\eta, \xi)_g:=\int_{\C^n}\eta\wedge \star \bar{\xi}=\int_{\C^n}\langle \eta, {\xi}\rangle_g \textup{vol}_g.\]
The differential operators $d:\Omega^{p-1}\left(\mathbb{C}^{n}\right)\rightarrow\Omega^{p}\left(\mathbb{C}^{n}\right)$, where $1\leq p\leq 2n$, have the following natural adjoint with respect to the bilinear form $(\cdot, \cdot)_g$.
\begin{defn}
The \textit{co-differential operator} $d^{*}:\Omega^{p}\left(\mathbb{C}^{n}\right)\rightarrow\Omega^{p-1}\left(\mathbb{C}^{n}\right)$
is defined by 
\begin{equation}
d^{*}\omega:=\left(-1\right)^{n\left(p+1\right)+1}\star d\star\omega,\label{eq:2.3}
\end{equation}
 where $\omega$ is any differential $p$-form. In particular, for $\C^2$ we have $d^{*}\omega=-\star d\star\omega$.
\end{defn}

\begin{defn}
A differential $p$-form $\omega$ is said to be \textit{Hodge-Laplace harmonic}
(HL-harmonic for short) if 
\begin{equation}
\Delta\omega:=\left(dd^{*}+d^{*}d\right)\omega=0,\label{eq:2.4}
\end{equation}
where $\Delta$ is referred to as the \textit{Hodge Laplacian} (or the Laplace-de Rham operator).
\end{defn}
For more information on complex differential forms we refer readers to \cite{baez1994gauge, darling1994differential, holm2008geometric, tu2017differential}.

\section{Maxwell's Equations}

Maxwell equations have several equivalent formulations (\cite{baez1994gauge, fleisch2008student, maxwell1865viii, swanson2014path}). 
From the viewpoint of physics,
there are versions of Maxwell's equations based on electric and magnetic
potentials that allow one to solve the equations as a boundary value
problem within the realms of classical physics and quantum mechanics.
Quantum mechanics is not in the scope of this paper, so we shall
only consider Maxwell's equations in the classical sense.
Moreover, the space-time formulations of Maxwell's equations are primarily
used in high-energy physics and gravitational physics, in conjunction
with Einstein's theories of special relativity and general relativity.
From the viewpoint of mathematics, Maxwell's equations are important
in vector calculus, potential field theory, gauge theory, differential geometry, topology, and many other
studies (\cite{baez1994gauge, garrity2015electricity, holm2008geometric}). 

\subsection{Classical Version}

Let $\mathbf{E}=\left(E_{1},E_{2},E_{3}\right)$ and $\mathbf{B}=\left(B_{1},B_{2},B_{3}\right)$
be electric and magnetic fields, respectively, in a convex region
of $\mathbb{R}^{4}$ with Lorentzian signature $+---$\footnote{This is essentially the Minkowski space-time $\mathbb{R}^{1,3}$.},
where the $E_{j},B_{j}$ are scalar-valued functions of time and space.
In this vector space, we represent points as coordinate column
vectors $\left(x_{0},x_{1},x_{2},x_{3}\right)$, where $x_{0}=ct$ ($c$ being the speed of light) 
is the time variable and $x_{1},x_{2},x_{3}$ are the space variables $\footnote{For simplicity, we often set $c=1$.}$.
Maxwell's equations are given by the following partial differential
equations: \begin{subequations}

\begin{eqnarray}
\nabla\cdotp\mathbf{B} & = & 0,\label{eq:3.1a}\\
\nabla\times\mathbf{E}+\frac{\partial\mathbf{B}}{\partial t} & = & 0,\label{eq:3.1b}\\
\nabla\cdotp\mathbf{E} & = & \rho,\label{eq:3.1c}\\
\nabla\times\mathbf{B}-\frac{\partial\mathbf{E}}{\partial t} & = & \mathbf{J},\label{eq:3.1d}
\end{eqnarray}

\end{subequations} 

\noindent where the scalar $\rho$ is the electric charge density
and the vector $\mathbf{J}$ is the electric current density vector.
The equations (\ref{eq:3.1a}) and (\ref{eq:3.1b}) are homogeneous, while
the equations (\ref{eq:3.1c}) and (\ref{eq:3.1d}) are inhomogeneous. From vector calculus, $\nabla\cdotp\left(\nabla\times\mathbf{A}\right)=0$
for any smooth vector field $\mathbf{A}=(A_1, A_2, A_3)$ and $\nabla\times\left(\nabla\phi\right)=0$
for any scalar function $\phi$ \footnote{From physics viewpoint, a scalar field.}. 
Now by Poincar\'{e}'s lemma, in $\mathbb{R}^{3}$, if $\nabla\cdotp\mathbf{B}=0$,
then $\mathbf{B}=\nabla\times\mathbf{A}$ for some vector field $\mathbf{A}$.
This leads us to the potential field theory aspect of electromagnetism.
So in this context, one considers the magnetic vector potential $\mathbf{A}$
and the electric scalar potential $\phi$ such that 
\begin{align}
\mathbf{B} & =\nabla\times\mathbf{A},\label{eq:3.9}\\
\mathbf{E} & =-\frac{\partial\mathbf{A}}{\partial t}-\nabla\phi.\label{eq:3.10}
\end{align}
 Putting (\ref{eq:3.9}) and (\ref{eq:3.10}) in equations (\ref{eq:3.1c})
and (\ref{eq:3.1d}), we obtain 
\begin{equation}
\rho=-\frac{\partial}{\partial t}\left(\nabla\cdotp\mathbf{A}\right)-\Delta\phi,\label{eq:3.11}
\end{equation}
 where $\Delta$ in (\ref{eq:3.11}) is the {\em Laplacian} in $\mathbb{R}^{3}$. 

So if we let $\mathbf{A}^{'}=\mathbf{A}-\nabla\psi, \phi^{'}=\phi+\frac{\partial\psi}{\partial t}$, for any scalar field $\psi$,  then the 4-vector $\left(\phi^{'},\mathbf{A}^{'}\right)$ solves
(\ref{eq:3.9}) and (\ref{eq:3.10}), but not uniquely. In particular,
this gives one the gauge freedom to choose $\phi$ such that $\left(\phi,\mathbf{A}\right)$
satisfies the {\em Lorenz}\footnote{Not to be confused with Lorentz!} {\em gauge}:
\begin{equation}
\frac{\partial\phi}{\partial t}+\nabla\cdotp\mathbf{A}=0, \label{eq:3.12}
\end{equation}
which implies that 
\begin{align}
\rho & =\left(\frac{\partial^{2}}{\partial t^{2}}-\Delta\right)\phi,\label{eq:3.13}\\
\mathbf{J} & =\left(\frac{\partial^{2}}{\partial t^{2}}-\Delta\right)\mathbf{A},\label{eq:3.14}
\end{align}
 where $\frac{\partial^{2}}{\partial t^{2}}-\Delta$ is the {\em d'Alembertian},
or {\em wave operator}. Since Lorentz transformations keep the Minkowski
metric invariant, the d'Alembertian gives a Lorentz scalar. Further, Maxwell's equations are Lorentz invariant and gauge
invariant.

\subsection{Differential Forms Version}

Set $\omega=\phi dx_{0}-A_{1}dx_{1}-A_{2}dx_{2}-A_{3}dx_{3},$ where ${\bf A}=(A_1, A_2, A_3)$ is a time-dependent smooth vector field in ${\mathbb R}^3$. This $1$-form $\omega$ is often
referred to as the {\em magnetic potential} 1-form. With the
Lorenz gauge from (\ref{eq:3.12}), we assume $\omega$ satisfies the following normalization:
\begin{equation}
\frac{\partial\phi}{\partial x_{0}}+\frac{\partial A_{1}}{\partial x_{1}}+\frac{\partial A_{2}}{\partial x_{2}}+\frac{\partial A_{3}}{\partial x_{3}}=0.\label{eq:3.15}
\end{equation}
 Now set 
\begin{equation}
J=\rho dx_{0}+J_{1}dx_{1}+J_{2}dx_{2}+J_{3}dx_{3}.\label{eq:3.16}
\end{equation}
Define the $2$-form $F_\omega=d\omega$. This is called the \textit{Faraday field strength} or simply the \textit{Faraday} $2$-\textit{form}.
Using (\ref{eq:3.9}) and (\ref{eq:3.10}), we have 
\begin{align}
F_\omega & =-dx_{0}\wedge\left(E_{1}dx_{1}+E_{2}dx_{2}+E_{3}dx_{3}\right)\label{eq:3.17}\\
 & -B_{1}dx_{2}\wedge dx_{3}+B_{2}dx_{1}\wedge dx_{3}-B_{3}dx_{1}\wedge dx_{2}.\nonumber 
\end{align}
 With the exterior derivative $d$ and the Hodge star operator $\star$, the Maxwell's equations take the concise form
\begin{align}
dF_\omega & =0,\label{eq:3.18}\\
\star \text{ } d \star F_\omega & =J,\label{eq:3.19}
\end{align}
where (\ref{eq:3.18}) is referred to as the \textit{Bianchi identity}. Equation (\ref{eq:3.18}) is equivalent to the homogeneous Maxwell's equations (\ref{eq:3.1a}) and (\ref{eq:3.1b}); while (\ref{eq:3.19}) is equivalent to the inhomogeneous Maxwell's equations (\ref{eq:3.1c}) and (\ref{eq:3.1d}). We may also refer to (\ref{eq:3.18}) 
and (\ref{eq:3.19}) as the exterior differential form of Maxwell's equations. And we say $\omega$, or alternatively $F_\omega$, is a solution to the Maxwell's equations if (\ref{eq:3.18}) and (\ref{eq:3.19}) are satisfied.

It is now apparent that solutions to the Maxwell's equations are not unique, since if $\omega$ is a solution then ${\omega'}$ is also a solution for $\omega'=\omega+d\psi$, where $\psi$ is any smooth function. This fact gives the gauge freedom of choosing $\phi$ that satisfies the Lorenz gauge condition (\ref{eq:3.12}). This differential form formulation of Maxwell's equations also makes the following remark apparent.

\begin{rem}
In vacuum, that is, when $\rho=0,J=0$, it follows from \textup{(\ref{eq:3.18})}
and \textup{(\ref{eq:3.19})} that every self-dual or anti-self-dual 2-form 
$F_\omega$ is a solution to the source-free Maxwell's equations.
\end{rem}

\section{Harmonic Solutions To Maxwell's Equations}

In this section we would like to consider Maxwell's equations in vacuum, 
and study some complex solutions and, in particular, harmonic solutions
to those equations. But before doing this, we shall first construct
a complex differential form solution to Maxwell's equations initially
not in vacuum. So let's begin by considering smooth functions $f_{1},f_{2},f_{\bar{1}},f_{\bar{2}}$
in the two variables $z_{1},z_{2}$. Define the complex differential form
\begin{equation}
\omega\left(z\right)=f_{1}dz_{1}+f_{2}dz_{2}+f_{\bar{1}}d\bar{z}_{1}+f_{\bar{2}}d\bar{z}_{2}.\label{eq:3.20}
\end{equation}
The form $\omega$ acts as a potential, much like the magnetic potential discussed in Section 3.2. In a more general context, it is often referred to as a {\em connection form}. The associated {\em curvature form} or {\em curvature field} is defined as $F_{\omega}:= d\omega + \omega\wedge \omega$. Since $\omega$ is scalar-valued, we have $\omega\wedge \omega=0$ and therefore $F_{\omega}=d\omega$. Direct computation shall verify that
\begin{align}
F_{\omega} & =\left(\partial_{1}f_{2}-\partial_{2}f_{1}\right)dz_{1}\wedge dz_{2}+\left(\bar{\partial}_{1}f_{\bar{2}}-\bar{\partial}_{2}f_{\bar{1}}\right)d\bar{z}_{1}\wedge d\bar{z}_{2}\nonumber \\
 & +\left(\partial_{1}f_{\bar{2}}-\bar{\partial}_{2}f_{1}\right)dz_{1}\wedge d\bar{z}_{2}+\left(\partial_{2}f_{\bar{1}}-\bar{\partial}_{1}f_{2}\right)dz_{2}\wedge d\bar{z}_{1}\label{eq:3.21}\\
 & +\left(\partial_{1}f_{\bar{1}}-\bar{\partial}_{1}f_{1}\right)dz_{1}\wedge d\bar{z}_{1}+\left(\partial_{2}f_{\bar{2}}-\bar{\partial}_{2}f_{2}\right)dz_{2}\wedge d\bar{z}_{2}.\nonumber 
\end{align}

Now in terms of real 2-forms, we can write (\ref{eq:3.21}) in the
following way:
\begin{align}
F_{\omega} & =-2i\left(\partial_{1}f_{\bar{1}}-\bar{\partial}_{1}f_{1}\right)dx_{0}\wedge dx_{1}\nonumber \\
 & +(\left(\partial_{1}f_{2}-\partial_{2}f_{1}\right)+\left(\bar{\partial}_{1}f_{\bar{2}}-\bar{\partial}_{2}f_{\bar{1}}\right)\nonumber \\
 & +\left(\partial_{1}f_{\bar{2}}-\bar{\partial}_{2}f_{1}\right)-\left(\partial_{2}f_{\bar{1}}-\bar{\partial}_{1}f_{2}\right))dx_{0}\wedge dx_{2}\nonumber \\
 & +i(\left(\partial_{1}f_{2}-\partial_{2}f_{1}\right)-\left(\bar{\partial}_{1}f_{\bar{2}}-\bar{\partial}_{2}f_{\bar{1}}\right)\nonumber \\
 & -\left(\partial_{1}f_{\bar{2}}-\bar{\partial}_{2}f_{1}\right)-\left(\partial_{2}f_{\bar{1}}-\bar{\partial}_{1}f_{2}\right))dx_{0}\wedge dx_{3}\label{eq:3.26}\\
 & -2i\left(\partial_{2}f_{\bar{2}}-\bar{\partial}_{2}f_{2}\right)dx_{2}\wedge dx_{3}\nonumber \\
 & +(-\left(\partial_{1}f_{2}-\partial_{2}f_{1}\right)-\left(\bar{\partial}_{1}f_{\bar{2}}-\bar{\partial}_{2}f_{\bar{1}}\right)\nonumber \\
 & +\left(\partial_{1}f_{\bar{2}}-\bar{\partial}_{2}f_{1}\right)-\left(\partial_{2}f_{\bar{1}}-\bar{\partial}_{1}f_{2}\right))dx_{1}\wedge dx_{3}\nonumber \\
 & +i(\left(\partial_{1}f_{2}-\partial_{2}f_{1}\right)-\left(\bar{\partial}_{1}f_{\bar{2}}-\bar{\partial}_{2}f_{\bar{1}}\right)\nonumber \\
 & +\left(\partial_{1}f_{\bar{2}}-\bar{\partial}_{2}f_{1}\right)+\left(\partial_{2}f_{\bar{1}}-\bar{\partial}_{1}f_{2}\right))dx_{1}\wedge dx_{2}.\nonumber 
\end{align}
One can easily determine the electric and magnetic components in view of (\ref{eq:3.17}): 
\begin{align*}
E_{1} & =2i\left(\partial_{1}f_{\bar{1}}-\bar{\partial}_{1}f_{1}\right),\\
E_{2} & =-\left(\left(\partial_{1}f_{2}-\partial_{2}f_{1}\right)+\left(\bar{\partial}_{1}f_{\bar{2}}-\bar{\partial}_{2}f_{\bar{1}}\right)+\left(\partial_{1}f_{\bar{2}}-\bar{\partial}_{2}f_{1}\right)-\left(\partial_{2}f_{\bar{1}}-\bar{\partial}_{1}f_{2}\right)\right),\\
E_{3} & =-i\left(\left(\partial_{1}f_{2}-\partial_{2}f_{1}\right)-\left(\bar{\partial}_{1}f_{\bar{2}}-\bar{\partial}_{2}f_{\bar{1}}\right)-\left(\partial_{1}\bar{f}_{2}-\bar{\partial}_{2}f_{1}\right)-\left(\partial_{2}\bar{f}_{1}-\bar{\partial}_{1}f_{2}\right)\right),\\
B_{1} & =2i\left(\partial_{2}f_{\bar{2}}-\bar{\partial}_{2}f_{2}\right),\\
B_{2} & =-\left(\partial_{1}f_{2}-\partial_{2}f_{1}\right)-\left(\bar{\partial}_{1}f_{\bar{2}}-\bar{\partial}_{2}f_{\bar{1}}\right)+\left(\partial_{1}f_{\bar{2}}-\bar{\partial}_{2}f_{1}\right)-\left(\partial_{2}f_{\bar{1}}-\bar{\partial}_{1}f_{2}\right),\\
B_{3} & =-i\left(\left(\partial_{1}f_{2}-\partial_{2}f_{1}\right)-\left(\bar{\partial}_{1}f_{\bar{2}}-\bar{\partial}_{2}f_{\bar{1}}\right)+\left(\partial_{1}f_{\bar{2}}-\bar{\partial}_{2}f_{1}\right)+\left(\partial_{2}f_{\bar{1}}-\bar{\partial}_{1}f_{2}\right)\right).
\end{align*}

 In general, since $ {\bf E}$ and ${\bf B}$ are complex vectors in $\C^3$, the ``electromagnetic dynamics'' occurs in $7$-dimensional space ($\dim \C^3$ plus one dimension for time). 
The spatial inner product is defined by
\[\langle {\bf E}, {\bf B}\rangle=E_1\bar{B_1}+E_2\bar{B_2}+E_3\bar{B_3}.\]
To facilitate the computation of this inner product and the energy density $\frac{1}{2}(|{\bf E}|^2+|{\bf B}|^2)$, we write (\ref{eq:3.21}) as 
 \begin{equation}\label{ce}
 F_\omega=F_{12}dz_1\wedge dz_2+F_{\bar{1}\bar{2}}d\bar{z}_{1}\wedge d\bar{z}_{2}+\sum_{j,k=1,2}F_{j\bar{k}}dz_j\wedge d{\overline{z_k}}.\end{equation}
 
   Then the following can be verified by direct computation:
\begin{align}
\langle {\bf E}, {\bf B}\rangle = 4F_{1\bar{1}}\overline{F _{2\bar{2}}}+2\left[(F_{12}-F_{2\bar{1}})(\overline{F_{12}+F_{2\bar{1}}})
+(F_{\bar{1}\bar{2}}+F_{1\bar{2}})(\overline{F_{\bar{1}\bar{2}}-F_{1\bar{2}}})\right].\label{eq:angle}
\end{align}
The energy density of the electromagnetic dynamics can be computed as 
\begin{align}
\frac{1}{2}(|{\bf E}|^2+|{\bf B}|^2)=2(|F_{12}|^2+|F_{\bar{1}\bar{2}}|^2 + \sum_{j,k=1,2}|F_{j\bar{k}}|^2).\label{eq:energy}
\end{align}
If $g$ is the Euclidean metric on $\C^2$, then one can also verify that \[|{\bf E}|^2+|{\bf B}|^2=\langle F_\omega, F_\omega\rangle.\]

In particular, one observes that if $f_{\bar{i}}=\bar{f}_i, i=1, 2$ then both the electric field ${\bf E}$ and the magnetic field ${\bf B}$ are real vector-fields in ${\mathbb R}^3$. For example, one can write \[E_1=2iF_{1\bar{1}}=2i\left(\overline{\bar{\partial}_{1}f_1}-\bar{\partial}_{1}f_{1}\right)\] and see that $E_1$ is $4$ times the imaginary part of $\bar{\partial}_{1}f_{1}$. Other components can be checked similarly. Alternatively, one may observe that since $\omega$ is real in this case, the curvature field $F_\omega=d\omega$ must be real. 
 In this case, one has $F_{\bar{1}\bar{2}}=\overline{F_{12}}$ and $F_{1\bar{2}}=-\overline{F_{2\bar{1}}}$, and it follows that
\begin{align}
\langle {\bf E}, {\bf B}\rangle = 4\left(F_{1\bar{1}}\overline{F _{2\bar{2}}}+|F_{12}|^2-|F_{2\bar{1}}|^2\right).\label{rangle}
\end{align}
If ${\bf E}$ and ${\bf B}$ are real vectors, then $\cos^{-1}\frac{ \langle {\bf E}, {\bf B}\rangle}{|{\bf E}||{\bf B}|}$ is the angle between the electric field and the magnetic field. But when ${\bf E}$ and ${\bf B}$ are complex vectors, the physical meaning of 
$\langle {\bf E}, {\bf B}\rangle$ is less clear.

Although the original Maxwell's equations were formulated under the Minkowski metric, the differential form formulation 
\textup{(\ref{eq:3.18})} and \textup{(\ref{eq:3.19})} makes good sense under other metrics. In the sequel we shall consider the Maxwell's equations in this formulation under two different metrics.

\subsection{Euclidean Metric Case}
Under the Euclidean metric, applying the Hodge star operator to (\ref{eq:3.21}) (cf. Example \ref{duals}), one has
\begin{align}
\star F_{\omega} & =\left(\partial_{1}f_{2}-\partial_{2}f_{1}\right)dz_{1}\wedge dz_{2}+\left(\bar{\partial}_{1}f_{\bar{2}}-\bar{\partial}_{2}f_{\bar{1}}\right)d\bar{z}_{1}\wedge d\bar{z}_{2}\nonumber \\
 & -\left(\partial_{1}f_{\bar{2}}-\bar{\partial}_{2}f_{1}\right)dz_{1}\wedge d\bar{z}_{2}-\left(\partial_{2}f_{\bar{1}}-\bar{\partial}_{1}f_{2}\right)dz_{2}\wedge d\bar{z}_{1}\label{eq:3.22}\\
 & +\left(\partial_{2}f_{\bar{2}}-\bar{\partial}_{2}f_{2}\right)dz_{1}\wedge d\bar{z}_{1}+\left(\partial_{1}f_{\bar{1}}-\bar{\partial}_{1}f_{1}\right)dz_{2}\wedge d\bar{z}_{2}.\nonumber 
\end{align}
It is immediate that
$F_{\omega}$ is self-dual if and only if  
\begin{align}
\partial_{1}f_{\bar{1}}-\bar{\partial}_{1}f_{1} &= \partial_{2}f_{\bar{2}}-\bar{\partial}_{2}f_{2}, \label{eq:3.26a}\\
\partial_{1}f_{\bar{2}}-\bar{\partial}_{2}f_{1} &= \partial_{2}f_{\bar{1}}-\bar{\partial}_{1}f_{2}=0,  \nonumber
\end{align}
and it is anti-self-dual if and only if 
\begin{align}
\partial_{1}f_{\bar{1}}-\bar{\partial}_{1}f_{1} &= -\left(\partial_{2}f_{\bar{2}}-\bar{\partial}_{2}f_{2}\right), \label{eq:3.26b}\\
\partial_{1}f_{2}-\partial_{2}f_{1} &=\bar{\partial}_{1}f_{\bar{2}}-\bar{\partial}_{2}f_{\bar{1}}=0. \nonumber
\end{align}
The following fact follows immediately from (\ref{eq:angle}).
\begin{cor}
Let $F_\omega$ be a self-dual solution to the Maxwell equations in vaccum with respect to the Euclidean metric on $\C^2$. Then $\langle {\bf E}, {\bf B}\rangle\geq 0$ holds with equality only if $\omega$ is a trivial solution.
\end{cor}
\proof In fact, Equations (\ref{eq:3.26a}) shows that $F_{1\bar{1}}=F_{2\bar{2}}$ and $F_{1\bar{2}}=F_{2\bar{1}}=0$. Hence by (\ref{eq:angle}) one has 
\[\langle {\bf E}, {\bf B}\rangle=2(2|F_{1\bar{1}}|^2+|F_{12}|^2+|F_{\bar{1}\bar{2}}|^2)\geq 0.\]
If $\langle {\bf E}, {\bf B}\rangle=0$ then all six coefficients of the $2$-forms in $F_\omega$ are $0$ and hence $\omega$ is a trivial solution.\qed\\

Since $F_{\omega}=d\omega$, $F_{\omega}$ is exact. Hence if $F_{\omega}$ is self-dual or anti-self-dual then \[d\star F_{\omega}= \pm d F_{\omega}=\pm d^{2}\omega=0.\] Now if $f_{1},f_{2}$ are holomorphic and $f_{\bar{1}}, f_{\bar{2}}$ are conjugate holomorphic, then the equations in (\ref{eq:3.26a}) are automatically satisfied. We thus have the following fact.
\begin{prop}\label{hol1}
The form $F_{\omega}$, determined by $\omega$ in \textup{(\ref{eq:3.20})} with $f_{1}, f_{2}$ holomorphic and $f_{\bar{1}}, f_{\bar{2}}$ conjugate holomorphic, is a self-dual solution for the source-free Maxwell's equations \textup{(\ref{eq:3.21})} and \textup{(\ref{eq:3.22})} under the Euclidean metric. Further, in this case
\begin{align}
F_{\omega}&=\left(\partial_{1}f_{2}-\partial_{2}f_{1}\right)dz_{1}\wedge dz_{2}+\left(\bar{\partial}_{1}f_{\bar{2}}-\bar{\partial}_{2}f_{\bar{1}}\right)d\bar{z}_{1}\wedge d\bar{z}_{2}.\label{eq:3.27}
\end{align}
\end{prop}
Consequently, we have
\begin{cor}\label{hol2}
If $\omega$ is of the form \textup{(\ref{eq:3.20})} and $f_{1}, f_{2}$ are holomorphic and $f_{\bar{1}}, f_{\bar{2}}$ are conjugate holomorphic, then
\begin{align*}
E_{1} & =B_{1}=0,\\
E_{2} & =B_{2}=-\left((\partial_{1}f_{2}-\partial_{2}f_{1})+(\bar{\partial}_{1}f_{\bar{2}}-\bar{\partial}_{2}f_{\bar{1}})\right),\\
E_{3} & =B_{3}=-i\left((\partial_{1}f_{2}-\partial_{2}f_{1})-(\bar{\partial}_{1}f_{\bar{2}}-\bar{\partial}_{2}f_{\bar{1}})\right).
\end{align*}
\end{cor}
It is surprising that in this case the electric field and the magnetic field coincide, or in other words they are mathematically indistinguishable. The same phenomenon occurs later in Example \ref{monople} on Dirac monople.

Of course self-duality (resp. anti-self-duality) of the differential form $F_{\omega}$ does not necessarily require $f_{j},f_{\bar{j}},\text{ }j=1,2$ being holomorphic (resp. conjugate holomorphic).
For instance, $f_{1}=z_{1}+\bar{z}_{1}, f_{2}=z_{2}+\bar{z}_{2}, f_{\bar{1}}=\bar{z}_{1}-z_{1}, f_{\bar{2}}=\bar{z}_{2}-z_{2}$ satisfy (\ref{eq:3.26a}), while
$f_{1}=z_{1}-\bar{z}_{2}, f_{2}=z_{2}-\bar{z}_{1}, f_{\bar{1}}=\bar{z}_{1}-z_{1}, f_{\bar{2}}=z_{2}-\bar{z}_{2}$ satisfy (\ref{eq:3.26b}). The next example exihibits a $1$-dimensional electromagnetic dynamics.

\begin{example}
Let $\tau\left(z\right)=f_{1}\left(dz_{1}+ d\bar{z}_{1}\right)+f_{2}\left(dz_{2}+ d\bar{z}_{2}\right)$, where $f_{1}$ and $f_{2}$ are holomorphic in $z_{1}$ and $z_{2}$. We claim that ${F}:=d\tau$ is self-dual if only if ${F}=m\left(dz_{1}\wedge d\bar{z}_{1}+dz_{2}\wedge d\bar{z}_{2}\right)$ for some constant $m$. First, it is clear that if ${F}=m\left(dz_{1}\wedge d\bar{z}_{1}+dz_{2}\wedge d\bar{z}_{2}\right)$
then it is self-dual. Conversely, since $f_{1},f_{2}$ are holomorphic, using the definition of ${F}$, after some simplification
we have 
\begin{align}   
{F} & =\partial_{1}f_{1}dz_{1}\wedge d\bar{z}_{1}+\left(\partial_{1}f_{2}-\partial_{2}f_{1}\right)dz_{1}\wedge dz_{2}+\partial_{1}f_{2}dz_{1}\wedge d\bar{z}_{2}\label{eq:3.28}\\
 & +\partial_{2}f_{1}dz_{2}\wedge d\bar{z}_{1}+\partial_{2}f_{2}dz_{2}\wedge d\bar{z}_{2}.\nonumber 
\end{align}
 After setting ${F}=\star {F}$ and comparing coefficients, it follows that 
\begin{equation}
\partial_{1}f_{1}=\partial_{2}f_{2},\text{ }\partial_{1}f_{2}=0=\partial_{2}f_{1}.\label{eq:3.29}
\end{equation}
The second relation in \textup{(\ref{eq:3.29})} implies $f_{1}$ is independent of $z_{2}$, and likewise $f_{2}$ is independent of $z_{1}$.
From the first relation, since $\partial_{1}f_{1}=\partial_{2}f_{2}$, $f_{1}$ and $f_{2}$ must be constant. Therefore, we can write $f_{1}=mz_{1}+c_{1}, f_{2}=mz_{2}+c_{2}$ for constants $m, c_{1}, c_{2}\in\mathbb{C}$.
 Substituting the conditions in \textup{(\ref{eq:3.29})} into \textup{(\ref{eq:3.28})},
we have 
\begin{align}
{F}&=m\left(dz_{1}\wedge d\bar{z}_{1}+dz_{2}\wedge d\bar{z}_{2}\right),\nonumber\\
&=-2im\left(dx_{0}\wedge dx_{1}+dx_{2}\wedge dx_{3}\right).\label{eq:3.30}
\end{align}
One can read off the electric and magnetic components in \textup{(\ref{eq:3.30})} by comparing with \textup{(\ref{eq:3.17})}:
\begin{align*}
E_{1} & =B_{1}=-2im,\\
E_{2} & =B_{2}=E_{3}=B_{3}=0,
\end{align*}
 which shows that the electromagnetic dynamics in this case is $1$-dimensional.
\end{example}

We now continue with the discussion on the general case that $f_{i}, f_{\bar{i}}, i=1, 2$ are smooth functions on $\C^2$. After applying the exterior derivative to (\ref{eq:3.22}) and some simplification we obtain
\begin{align}
d\star F_{\omega} & =\left(-\partial_{2}\bar{\partial}_{1}f_{1}+\left(2\partial_{1}\bar{\partial}_{1}+\partial_{2}\bar{\partial}_{2}\right)f_{2}-\partial_{1}\partial_{2}f_{\bar{1}}-\partial_{2}^{2}f_{\bar{2}}\right)dz_{1}\wedge dz_{2}\wedge d\bar{z}_{1}\nonumber \\
 & +\left(-\left(\partial_{1}\bar{\partial}_{1}+2\partial_{2}\bar{\partial}_{2}\right)f_{1}+\partial_{1}\bar{\partial}_{2}f_{2}+\partial_{1}^{2}f_{\bar{1}}+\partial_{1}\partial_{2}f_{\bar{2}}\right)dz_{1}\wedge dz_{2}\wedge d\bar{z}_{2}\label{eq:3.23}\\
 & +\left(-\bar{\partial}_{1}\bar{\partial}_{2}f_{1}-\bar{\partial}_{2}^{2}f_{2}-\bar{\partial}_{2}\partial_{1}f_{\bar{1}}+\left(2\bar{\partial}_{1}\partial_{1}+\bar{\partial}_{2}\partial_{2}\right)f_{\bar{2}}\right)dz_{1}\wedge d\bar{z}_{1}\wedge d\bar{z}_{2}\nonumber \\
 & +\left(\bar{\partial}_{1}^{2}f_{1}+\bar{\partial}_{1}\bar{\partial}_{2}f_{2}-\left(\bar{\partial}_{1}\partial_{1}+2\bar{\partial}_{2}\partial_{2}\right)f_{\bar{1}}+\bar{\partial}_{1}\partial_{2}f_{\bar{2}}\right)dz_{2}\wedge d\bar{z}_{1}\wedge d\bar{z}_{2}.\nonumber 
\end{align}
 
Applying the Hodge star operator to (\ref{eq:3.23}) and rearranging
terms, we have 
\begin{align}
\star \text{ } d\star F_{\omega} & =2\left(\left(\partial_{1}\bar{\partial}_{1}+2\partial_{2}\bar{\partial}_{2}\right)f_{1}-\partial_{1}\bar{\partial}_{2}f_{2}-\partial_{1}^{2}f_{\bar{1}}-\partial_{1}\partial_{2}f_{\bar{2}}\right)dz_{1}\nonumber \\
 & +2\left(-\bar{\partial}_{1}^{2}f_{1}-\bar{\partial}_{1}\bar{\partial}_{2}f_{2}+\left(\bar{\partial}_{1}\partial_{1}+2\bar{\partial}_{2}\partial_{2}\right)f_{\bar{1}}-\bar{\partial}_{1}\partial_{2}f_{\bar{2}}\right)d\bar{z}_{1}\label{eq:3.24}\\
 & +2\left(-\partial_{2}\bar{\partial}_{1}f_{1}+\left(2\partial_{1}\bar{\partial}_{1}+\partial_{2}\bar{\partial}_{2}\right)f_{2}-\partial_{1}\partial_{2}f_{\bar{1}}-\partial_{2}^{2}f_{\bar{2}}\right)dz_{2}\nonumber \\
 & +2\left(-\bar{\partial}_{1}\bar{\partial}_{2}f_{1}-\bar{\partial}_{2}^{2}f_{2}-\bar{\partial}_{2}\partial_{1}f_{\bar{1}}+\left(2\bar{\partial}_{1}\partial_{1}+\bar{\partial}_{2}\partial_{2}\right)f_{\bar{2}}\right)d\bar{z}_{2}.\nonumber 
\end{align}
The RHS of (\ref{eq:3.24}) can be viewed as a complex current
form that we may denote by $J$. For the 
sake of simplicity, we can rewrite (\ref{eq:3.24}) as 
\begin{equation}
\star \text{ } d\star F_{\omega}=J:=P_{1}dz_{1}+P_{\bar{1}}d\bar{z}_{1}+P_{2}dz_{2}+P_{\bar{2}}d\bar{z}_{2},\label{eq:3.25}
\end{equation}
 where the coefficients $P_{j}, P_{\bar{j}}, j=1, 2$ can be easily read off from (\ref{eq:3.24}).
 Likewise, we can write the RHS of (\ref{eq:3.24}) in terms of real
1-forms in view of (\ref{eq:3.16}): 
\[
J=\left(P_{1}+P_{\bar{1}}\right)dx_{0}+i\left(P_{1}-P_{\bar{1}}\right)dx_{1}+\left(P_{2}+P_{\bar{2}}\right)dx_{2}+i\left(P_{2}-P_{\bar{2}}\right)dx_{3}.
\]
 Again, if $f_{\bar{i}}=\bar{f}_i, i=1, 2$ then $J$ is real. Here $P_{1}+P_{\bar{1}}$ can be considered
as a scalar electric charge density $\rho$, while the last three
coefficients for the above current 1-form correspond to the last
three components of the electric current density vector $\mathbf{J}=\left(\rho,J_{1},J_{2},J_{3}\right)$
in Section 3.1. So to summarize, with smooth complex-valued
functions $f_{1},f_{2},f_{\bar{1}},f_{\bar{2}}$ in two variables $z_{1},z_{2}$, we determined
from the complex differential form $\omega$ in (\ref{eq:3.20}),
a solution $F_{\omega}$ (in (\ref{eq:3.21})) to a complex analogue of Maxwell's
equations $dF_{\omega}=0$ and $\star d\star F_{\omega}=J$.

Let $\nabla^{2}:= 4\left(\partial_{1}\bar{\partial}_{1}+ \partial_{2}\bar{\partial}_{2}\right)$ be the Laplacian on $\mathbb{C}^{2}$ in the Euclidean metric. Then a complex function $f$ is said to be $\textit{harmonic}$ if $\nabla^{2}f=0$ on $\mathbb{C}^{2}$. The following is the main result of this subsection.
\begin{thm}\label{thm:main1}
Let $f_{j},f_{\bar{j}},\text{ }j=1,2$ be harmonic functions. Then the complex differential form $\omega$ is a solution to the source-free Maxwell's equations in the Euclidean metric if and only if $\bar{\partial}_{1}f_{1}+\bar{\partial}_{2}f_{2}+\partial_{1}f_{\bar{1}}+\partial_{2}f_{\bar {2}}$ is constant.
\end{thm}
\proof
First,  we can rewrite (\ref{eq:3.23}) as
\begin{align}
d\star F_{\omega} &= \left(\nabla^{2}f_{2}-\partial_{2} \bar{\partial}_{1}f_{1}-\partial_{2}^{2}f_{\bar{2}}-\partial_{1}\partial_{2}f_{\bar{1}}+\partial_{1}\bar{\partial}_{1}f_{2}\right) dz_{1}\wedge dz_{2}\wedge d\bar{z}_{1} \nonumber \\
& +\left(-\nabla^{2}f_{1}+\partial_{1}\bar{\partial}_{2}f_{2}+\partial_{1}^{2}f_{\bar{1}}+\partial_{1}\partial_{2}f_{\bar{2}}-\partial_{2}\bar{\partial}_{2}f_{1}\right) dz_{1}\wedge dz_{2}\wedge d\bar{z}_{2} \label{eq:3.31} \\
& + \left(\nabla^{2}f_{\bar{2}}-\partial_{1}\bar{\partial}_{2}f_{\bar{1}}-\bar{\partial}_{2}^{2}f_{2}+\partial_{1}\bar{\partial}_{1}f_{\bar{2}}-\bar{\partial}_{1}\bar{\partial}_{2}f_{1}\right) dz_{1}\wedge d\bar{z}_{1}\wedge d\bar{z}_{2} \nonumber \\
& +\left(-\nabla^{2}f_{\bar{1}}+\partial_{2}\bar{\partial}_{1}f_{\bar{2}}+\bar{\partial}_{1}^{2}f_{1}-\partial_{2}\bar{\partial_{2}}f_{\bar{1}}+\bar{\partial}_{1}\bar{\partial}_{2}f_{2}\right) dz_{2}\wedge d\bar{z}_{1}\wedge d\bar{z}_{2}. \nonumber 
\end{align}
Since $f_{j},f_{\bar{j}},\text{ }j=1,2$ are harmonic, we have $\nabla^{2}f_{j}=0$
and $\nabla^{2}f_{\bar{j}}=0$, which implies
\begin{align*} 
d\star F_{\omega} &= \left(-\partial_{2} \bar{\partial}_{1}f_{1}-\partial_{2}^{2}f_{\bar{2}}-\partial_{1}\partial_{2}f_{\bar{1}}-\partial_{2}\bar{\partial}_{2}f_{2}\right) dz_{1}\wedge dz_{2}\wedge d\bar{z}_{1}\\
& +\left(\partial_{1}\bar{\partial}_{2}f_{2}+\partial_{1}^{2}f_{\bar{1}}+\partial_{1}\partial_{2}f_{\bar{2}}+\partial_{1}\bar{\partial}_{1}f_{1}\right) dz_{1}\wedge dz_{2}\wedge d\bar{z}_{2}\\
& + \left(-\partial_{1}\bar{\partial}_{2}f_{\bar{1}}-\bar{\partial}_{2}^{2}f_{2}-\partial_{2}\bar{\partial}_{2}f_{\bar{2}}-\bar{\partial}_{1}\bar{\partial_{2}}f_{1}\right)  dz_{1}\wedge d\bar{z}_{1}\wedge d\bar{z}_{2}\\
& +\left(\partial_{2}\bar{\partial}_{1}f_{\bar{2}}+\bar{\partial}_{1}^{2}f_{1}+\partial_{1}\bar{\partial}_{1}f_{\bar{1}}+\bar{\partial}_{1}\bar{\partial}_{2}f_{2}\right) dz_{2}\wedge d\bar{z}_{1}\wedge d\bar{z}_{2}.
\end{align*}

For $F_{\omega}$ to satisfy the source-free Maxwell's equations, we require $d\star F_{\omega}=0$, which we may now write in matrix form as
\[\text{diag}\{-\partial_{2},\partial_ {1},-\bar{\partial}_{2},\bar{\partial}_{1}\}
\begin{pmatrix}
\bar{\partial_{1}} & \bar{\partial_{2}} & \partial_{1} & \partial_{2}\\
\bar{\partial_{1}} & \bar{\partial_{2}} & \partial_{1} & \partial_{2}\\
\bar{\partial_{1}} & \bar{\partial_{2}} & \partial_{1} & \partial_{2}\\
\bar{\partial_{1}} & \bar{\partial_{2}} & \partial_{1} & \partial_{2}
\end{pmatrix}
\begin{pmatrix} f_{1} \\ f_{2} \\ f_{\bar{1}} \\ f_{\bar{2}}
\end{pmatrix} = \mathbf{0},
\]
where ``diag'' stands for a diagonal matrix and $\mathbf{0}$ denotes the column 4-vector of zeroes. Clearly, the above is true if and only if $\bar{\partial}_{1}f_{1}+\bar{\partial}_{2}f_{2}+\partial_{1}f_{\bar{1}}+\partial_{2}f_{\bar {2}}$ is a constant.\qed \\

\begin{example}\label{monople}
With $z_{1}=x_{0}+ix_{1}, z_{2}=x_{2}+ix_{3}$, consider the form
\begin{align}
\omega\left(z\right) &= i\left(x_{0}dx_{1}-x_{1}dx_{0}+x_{2}dx_{3}-x_{3}dx_{2}\right) \nonumber\\
 &= \frac{1}{2}\left(\eta\left(z\right)- \overline{\eta\left(z\right)}\right), \label{eq:3.32}
\end{align}
where $\eta\left(z\right)=\bar{z}_{1}dz_{1}+\bar{z}_{2}dz_{2}$. In this case, we have $f_{1}=\frac{1}{2}\bar{z}_{1}, f_{2}=\frac{1}{2}\bar{z}_{2}, f_{\bar{1}}=-\frac{1}{2}z_{1}, f_{\bar{2}}=-\frac{1}{2}z_{2}$, where $f_{j},f_{\bar{j}},\text{ }j=1,2$ are not holomorphic (resp. conjugate holomorphic) but are all clearly harmonic.  Moreover, we have the associated curvature form
\begin{equation*}
F_{\omega}= d\omega= -\left(dz_{1}\wedge d\bar{z}_{1}+dz_{2}\wedge d\bar{z}_{2}\right),
\end{equation*}
which is self-dual in the Euclidean metric. Now observe that							
\begin{align*}
\partial_{1}f_{\bar{1}}-\bar{\partial}_{1}f_{1} &= -1 = \partial_{2}f_{\bar{2}}-\bar{\partial}_{2}f_{2}, \\
\partial_{1}f_{\bar{2}} &= 0 = \bar{\partial}_{1}f_{2},\\
\partial_{2}f_{\bar{1}} &= 0 = \bar{\partial}_{2}f_{1}.
\end{align*}
The above conditions clearly satisfy the self-duality of $F_{\omega}$. Moreover, $\bar{\partial}_{1}f_{1}+\bar{\partial}_{2}f_{2}+\partial_{1}f_{\bar{1}}+\partial_{2}f_{\bar {2}}=0$. 
\end{example}

The above example is related to the {\em Dirac monopole} (\cite{ritter2003gauge}), which is a hypothetical magnetic charge. The original idea was proposed in a 1931 paper by Paul Dirac \textup{(\cite{dirac1931quantised})}. Evidently, the above example shows that the existence of Dirac monopoles does not conflict with Maxwell's equations in vacuum (away from the magnetic monople). See \cite{dirac1931quantised, pinfold2010dirac} and the references therein for more background on this rather intriguing subject. A notable fact here is that ${\bf E}={\bf B}$. It is easy to compute that ${\bf E}={\bf B}=(-2i, 0, 0)$, and therefore
\[\langle {\bf E}, {\bf B}\rangle =\frac{1}{2}(|{\bf E}|^2+|{\bf B}|^2)=4.\]

Theorem \ref{thm:main1} leads to easy constructions of non-self-dual solutions to Equations \textup{(\ref{eq:3.18})} and \textup{(\ref{eq:3.19})} in vaccum.
\begin{example}
In fact, a simple working example that satisfies the conditions in Theorem \ref{thm:main1}, but fails the conditions for $F_{\omega}$ to be  self-dual nor anti-self-dual, is $f_{1}=2\bar{z}_{1}-z_{2}, f_{2}=z_{1}+2\bar{z}_{2}, f_{\bar{1}}=z_{1}+\bar{z}_{1}, f_{\bar{2}}=z_{2}+\bar{z}_{2}$.
\end{example}

\subsection{Minkowski Metric Case}

We believe that much of the work in Section 4.1 can be done in a parallel manner with respect to other bilinear forms $\langle \cdot, \cdot\rangle_g$ on ${\mathbb R}^4$, where $g$ is a nondegenerate constant $4\times 4$ self-adjoint matrix. But since 
the Minkowski metric on ${\mathbb R}^{1,3}$ is more conforming with our reality, and it is indeed where the Maxwell's equations were initially studied, we shall work it out in details in this subsection. 

Recall that the d'Alembertian in the Minkowski metric is given by
\begin{equation}
\Box = \frac{\partial^{2}}{\partial x_{0}^{2}}- \frac{\partial^{2}}{\partial x_{1}^{2}}- \frac{\partial^{2}}{\partial x_{2}^{2}}- \frac{\partial^{2}}{\partial x_{3}^{2}} = \frac{\partial^{2}}{\partial x_{0}^{2}}- \nabla^2,
\label{eq:3.19b}
\end{equation}
where again $\nabla^2$ is the Laplacian in $\mathbb{R}^{3}$. Using the identities
\begin{align}
&\partial_{1} = \frac{1}{2}\left(\frac{\partial}{\partial x_{0}}-i\frac{\partial}{\partial x_{1}}\right), \hspace{1cm}  \bar{\partial}_{1} = \frac{1}{2}\left(\frac{\partial}{\partial x_{0}}+i\frac{\partial}{\partial x_{1}}\right) \nonumber \\
&\partial_{2} = \frac{1}{2}\left(\frac{\partial}{\partial x_{2}}-i\frac{\partial}{\partial x_{3}}\right), \hspace{1cm} \bar{\partial}_{2} = \frac{1}{2}\left(\frac{\partial}{\partial x_{2}}+i\frac{\partial}{\partial x_{3}}\right),\label{eq:3.19c}
\end{align}
we have the following.
\begin{lem}
In terms of complex variables, the d'Alembertian in the Minkowski metric is given by
\begin{equation}
\Box = 2\left(\partial_{1}^{2}+ \bar{\partial}_{1}^{2}-2 \partial_{2}\bar{\partial}_{2}\right).
\label{eq:3.19d}
\end{equation}
\end{lem}
\proof
Using the first two identities in (\ref{eq:3.19c}), note that
\begin{align*}
\partial_{1}+\bar{\partial}_{1} &= \frac{\partial}{\partial x_{0}}, \hspace{1cm}
i\left(\partial_{1}-\bar{\partial}_{1}\right) = \frac{\partial}{\partial x_{1}}.
\end{align*}
Then the first two terms of the d'Alembertian in (\ref{eq:3.19b}) are
\begin{align*}
\left(\frac{\partial}{\partial x_{0}}+\frac{\partial}{\partial x_{1}}\right)\left(\frac{\partial}{\partial x_{0}}-\frac{\partial}{\partial x_{1}}\right) &= \left(\partial_{1}+\bar{\partial}_{1}+i\left(\partial_{1}-\bar{\partial}_{1}\right)\right)\left(\partial_{1}+\bar{\partial}_{1}-i\left(\partial_{1}-\bar{\partial}_{1}\right)\right) \\
 &= \left(\partial_{1}+\bar{\partial}_{1}\right)^{2}+\left(\partial_{1}-\bar{\partial}_{1}\right)^{2} \\
 &= 2\left(\partial_{1}^{2} +\bar{\partial}_{1}^{2}\right).
\end{align*}
Similarly, using the last two identities in (\ref{eq:3.19c}), it follows that
 the last two terms of the d'Alembertian in (\ref{eq:3.19b}) are $-4\partial_{2}\bar{\partial}_{2}$.
Therefore, $\Box = 2\left(\partial_{1}^{2}+ \bar{\partial}_{1}^{2}-2 \partial_{2}\bar{\partial}_{2}\right)$.
\qed \\

\begin{defn}
A function $f$ is said to be $M$-harmonic if $\Box f=0$.
\end{defn}
Here the ``$M$'' refers to the Minkowski metric. It is well-known that $M$-harmonic functions $\psi$ describe waves propagating in ${\mathbb R}^{1,3}$. A smooth $1$-form as defined in \textup{(\ref{eq:3.20})} is said to be {\em wavelike} if the functions $f_j, f_{\bar{j}}, j=1, 2$ are all $M$-harmonic. Likewise, the curvature field $F_\omega$ is said to be wavelike if its coefficient functions $F_{12}, F_{\bar{1}\bar{2}}$, and $F_{j\bar{k}}, j, k=1, 2$ in (\ref{ce}) are all $M$-harmonic. It is easy to check that if $\omega$ is wavelike then $F_\omega$ is also wavelike, for instance,
\[\Box F_{12}=\Box \left(\partial_{1}f_{2}-\partial_{2}f_{1}\right)=\partial_{1}\Box f_{2}-\partial_{2}\Box f_{1}=0,\]
and other coefficients are checked similarly. In particular, this implies that $\Box{E_j}=\Box{B_j}=0, 1\leq j\leq 3$, i.e., ${\bf E}$ and ${\bf B}$ are waves in the space $\C^3$. However, as we will see a bit later, there exists non-wavelike $\omega$ for which $F_\omega$ is wavelike.

Now let $f_{j},f_{\bar{j}},\text{ }j=1,2$ be complex smooth functions as before. Then $F_{\omega}$ is still the same as in (\ref{eq:3.21}). However in the Minkowski metric, by Example \ref{mink} we have
\begin{align}
\star F_{\omega} & =\left(\bar{\partial}_{1}f_{\bar{2}}-\bar{\partial}_{2}f_{\bar{1}}\right)dz_{1}\wedge d\bar{z}_{2}+\left(\partial_{2}f_{1}-\partial_{1}f_{2}\right)dz_{2}\wedge d\bar{z}_{1}\nonumber \\
 & +\left(\partial_{2}f_{\bar{1}}-\bar{\partial}_{1}f_{2}\right)dz_{1}\wedge dz_{2}+\left(\bar{\partial}_{2}f_{1}-\partial_{1}f_{\bar{2}}\right)d\bar{z}_{1}\wedge d\bar{z}_{2}\label{eq:3.33}\\
 & +\left(\partial_{2}f_{\bar{2}}-\bar{\partial}_{2}f_{2}\right)dz_{1}\wedge d\bar{z}_{1}+\left(\bar{\partial}_{1}f_{1}-\partial_{1}f_{\bar{1}}\right)dz_{2}\wedge d\bar{z}_{2}.\nonumber 
\end{align}

In the Minkowski metric on $\mathbb{C}^{2}$, the curvature form $F_{\omega}$ is self-dual (resp. anti-self-dual) provided $\star F_{\omega} =\pm iF_{\omega}$ (\cite{baez1994gauge, felsager2012geometry}). So self-duality of $F_{\omega}$ requires
\begin{align}
\partial_{2}f_{\bar{1}}-\bar{\partial}_{1}f_{2} &= i\left(\partial_{1}f_{2}-\partial_{2}f_{1}\right), \nonumber\\
\bar{\partial}_{1}f_{\bar{2}}-\bar{\partial}_{2}f_{\bar{1}} &= i\left(\partial_{1}f_{\bar{2}}-\bar{\partial}_{2}f_{1}\right), \label{eq:3.34}\\
\partial_{2}f_{\bar{2}}-\bar{\partial}_{2}f_{2} &= i\left(\partial_{1}f_{\bar{1}}-\bar{\partial}_{1}f_{1}\right). \nonumber
\end{align}
On the other hand, anti-self-duality of $F_{\omega}$ requires
\begin{align}
\partial_{2}f_{\bar{1}}-\bar{\partial}_{1}f_{2} &= -i\left(\partial_{1}f_{2}-\partial_{2}f_{1}\right), \nonumber\\
\bar{\partial}_{1}f_{\bar{2}}-\bar{\partial}_{2}f_{\bar{1}} &= -i\left(\partial_{1}f_{\bar{2}}-\bar{\partial}_{2}f_{1}\right), \label{eq:3.35}\\
\partial_{2}f_{\bar{2}}-\bar{\partial}_{2}f_{2} &= -i\left(\partial_{1}f_{\bar{1}}-\bar{\partial}_{1}f_{1}\right). \nonumber
\end{align}

The following fact is immediate. 

\begin{cor}
If $\omega$ is a self-dual or anti-self-dual solution to the Maxwell's equations in vaccum with respect to the Minkowski metric, then
$\langle {\bf E}, {\bf B}\rangle$ is either $0$ or purely imaginary. In particular, if $\omega$ is a real self-dual or anti-self-dual solution then $\langle {\bf E}, {\bf B}\rangle=0$.
\end{cor}
\proof If $\omega$ is self-dual, then (\ref{eq:3.35}) indicates that 
\[F_{2\bar{1}}=iF_{12},\ \ F_{\bar{1}\bar{2}}=iF_{1\bar{2}},\ \ F_{2\bar{2}}=iF_{1\bar{1}}.\]
Applying these relations to (\ref{eq:angle}), one has
 \begin{align*}
 \langle {\bf E}, {\bf B}\rangle=& -4i |F_{1\bar{1}}|^2+2(1-i)^2 |F_{12}|^2+2(1+i)^2 |F_{1\bar{2}}|^2\\
&=-4i\left(|F_{1\bar{1}}|^2+ |F_{12}|^2- |F_{1\bar{2}}|^2\right).
\end{align*}
In the case $\omega$ is anti-self-dual, parallel computations yield 
\[\langle {\bf E}, {\bf B}\rangle=4i\left(|F_{1\bar{1}}|^2+ |F_{12}|^2- |F_{1\bar{2}}|^2\right).\]
If $\omega$ is real then $ \langle {\bf E}, {\bf B}\rangle$ is real and therefore it must be equal to $0$.\qed\\

To proceed, as in the previous subsection we assume $F_{\omega}$ is neither self-dual nor anti-self-dual in the Minkowski metric. So for the source-free Maxwell's equations to be satisfied, we require $d\star F_{\omega}=0$. In this case, we have
\begin{align}
d\star F_{\omega} &= \left(\partial_{1}\partial_{2}f_{1}-\left(\partial_{1}^{2}+\bar{\partial}_{1}^{2}-\partial_{2}\bar{\partial}_{2}\right)f_{2}+\partial_{2}\bar{\partial}_{1}f_{\bar{1}}-\partial_{2}^{2}f_{\bar{2}}\right) dz_{1}\wedge dz_{2}\wedge d\bar{z}_{1} \nonumber\\
& +\left(\partial_{1}\bar{\partial}_{1}f_{1}-\bar{\partial}_{1}\bar{\partial}_{2}f_{2}-\left(\partial_{1}^{2}-2\partial_{2}\bar{\partial}_{2}\right)f_{\bar{1}}-\partial_{2}\bar{\partial}_{1}f_{\bar{2}}\right) dz_{1}\wedge dz_{2}\wedge d\bar{z}_{2} \label{eq:3.36}\\
& +\left(\partial_{1}\bar{\partial}_{2}f_{1}-\bar{\partial}_{2}^{2}f_{2}+\bar{\partial}_{1}\bar{\partial}_{2}f_{\bar{1}}-\left(\partial_{1}^{2}+\bar{\partial}_{1}^{2}-\partial_{2}\bar{\partial}_{2}\right)f_{\bar{2}}\right) dz_{1}\wedge d\bar{z}_{1}\wedge d\bar{z}_{2} \nonumber\\
& +\left(-\left(\bar{\partial}_{1}^{2}-2\partial_{2}\bar{\partial}_{2}\right)f_{1}-\partial_{1}\bar{\partial}_{2}f_{2}+\partial_{1}\bar{\partial}_{1}f_{\bar{1}}-\partial_{1}\partial_{2}f_{\bar{2}}\right) dz_{2}\wedge d\bar{z}_{1}\wedge d\bar{z}_{2}.\nonumber
\end{align}

Observe that for $f_{j},f_{\bar{j}},\text{ }j=1,2$ holomorphic and respectively conjugate holomorphic, equations (\ref{eq:3.34}) and (\ref{eq:3.35}) imply that all the coefficients of the $2$-forms in $F_\omega$ are $0$. Hence there is no nontrivial self-dual or anti-self-dual solution to the source-free Maxwell equations in this case. However, it follows from the above computation that
\[d\star F_{\omega}=\left(\partial_{1}\partial_{2}f_{1}-\partial_{1}^{2}f_2\right)dz_{1}\wedge dz_{2}\wedge d\bar{z}_{1} +
\left(\bar{\partial}_{1}\bar{\partial}_{2}f_{\bar{1}}-\bar{\partial}_{1}^{2}f_{\bar{2}}\right)dz_{1}\wedge d\bar{z}_{1}\wedge d\bar{z}_{2}.\]
Hence $d\star F_{\omega}=0$ if and only if 
\begin{equation}\label{hdual}
\partial_{1}\partial_{2}f_{1}-\partial_{1}^{2}f_2=\bar{\partial}_{1}\bar{\partial}_{2}f_{\bar{1}}-\bar{\partial}_{1}^{2}f_{\bar{2}}=0.
\end{equation}
Further, if $f_{\bar{j}}=\overline{f_j},\text{ }j=1,2$ then the above two equations are the same. One thus obtains the following fact.
\begin{prop}\label{hs}
Let $f_1$ and $f_2$ be holomorphic functions and $f_{\bar{j}}=\overline{f_j},\text{ }j=1,2$. Then $\omega$ is a solution to the Maxwell's equations in vacuum with respect to the Minkowski metric if and only if $\partial_{2}f_{1}-\partial_{1}f_2$ is independent of the variable $z_1$.
\end{prop}
Similar to Proposition \ref{hol1} and Corollary \ref{hol2}, in this case we have 
\begin{align*}
F_{\omega}&=\left(\partial_{1}f_{2}-\partial_{2}f_{1}\right)dz_{1}\wedge dz_{2}+\overline{ \left(\partial_{1}f_{2}-\partial_{2}f_{1}\right)} d\bar{z}_{1}\wedge d\bar{z}_{2},\\
&=2\Re\left((\partial_{1}f_{2}-\partial_{2}f_{1})dz_{1}\wedge dz_{2}\right),
\end{align*}
and consequently,
\begin{align*}
E_{1} & =B_{1}=0,\\
E_{2} & =B_{2}=2\Re (\partial_{1}f_{2}-\partial_{2}f_{1}),\\
E_{3} & =B_{3}=-2\Im (\partial_{2}f_{1}-\partial_{1}f_{2}),
\end{align*}
where $\Re(a)$ and $\Im(a)$ stand for the real, and respectively, imaginary part of a complex number $a$. Observe that ${\bf E}={\bf B}$ in this case, which resembles the Dirac monople example we examined earlier.

\begin{example}\label{f1f2}
There are plenty of holomorphic functions $f_1$ and $f_2$ that satisfy the condition in the above proposition. For instance, let \[f_1=z_1^2h(z_2)+g(z_2),\ f_2=\frac{z_1^3}{3}\partial_2h(z_2),\] where $g$ and $h$ are arbitrary one-variable entire functions. Then
$\partial_{2}f_{1}-\partial_{1}f_2=\partial_2g(z_2)$, which is independent of $z_1$. Further, since in this case
\[\Box f_1=2h(z_2),\ \ \  \Box f_2=2z_1\partial_2h(z_2),\] which can be nonzero, the $1$-form $\omega$ may not be wavelike. 
Further, it is easy to see that $\partial_2g(z_2)$ is M-harmonic and hence $F_\omega$ is wavelike. 
\end{example}
We state this observation as follows.
\begin{cor} \label{nonwavelike}
There are real analytic non-wavelike solutions to the Maxwell's equations in vacuum.
\end{cor} 
\begin{rem} \textup{Corollary \ref{hol2}} and the above observations also indicate that, under both the Euclidean metric and the Minkowski metric, the Maxwell's equations in vacuum have solutions in which the electric field and the magnetic field are mathematically indistinguishable. However, it is not clear if such solutions exist in nature.
\end{rem}

Now coming back to our familiar wavelike solutions we have the following fact. Its proof is similar to that of Theorem \ref{thm:main1}.
\begin{thm}\label{thm:main2}
Assume $\omega$ as in \textup{(\ref{eq:3.20})} is wavelike. Then it is a solution to the Maxwell's equations in vacuum under the Minkowski metric if and only if  $\partial_{1}f_{1}-\bar{\partial}_{2}f_{2}+\bar{\partial}_{1}f_{\bar{1}}-\partial_{2}f_{\bar {2}}$ is constant.
\end{thm}
\proof
 Let $\mathbf{0}$ denote the column 4-vector of zeroes.
Since $\omega$ is wavelike, we have $\Box f_{j}=\Box f_{\bar{j}}=0, j=1, 2$, which means
\begin{equation}
\left(\partial_{1}^{2}+\bar{\partial}_{1}^{2} -2\partial_{2}\bar{\partial}_{2}\right)f_{j}=0, \left(\partial_{1}^{2}+\bar{\partial}_{1}^{2} -2\partial_{2}\bar{\partial}_{2}\right)f_{\bar{j}}=0.
\label{eq:3.37}
\end{equation}
Plugging (\ref{eq:3.37}) into (\ref{eq:3.36}) and rearranging terms, we have
\begin{align*}
d\star F_{\omega} &= \left(\partial_{1}\partial_{2}f_{1}-\partial_{2}\bar{\partial}_{2}f_{2}+\partial_{2}\bar{\partial}_{1}f_{\bar{1}}-\partial_{2}^{2}f_{\bar{2}}\right) dz_{1}\wedge dz_{2}\wedge d\bar{z}_{1} \\
& +\left(\partial_{1}\bar{\partial}_{1}f_{1}-\bar{\partial}_{1}\bar{\partial}_{2}f_{2}+\bar{\partial}_{1}^{2}f_{\bar{1}}-\partial_{2}\bar{\partial}_{1}f_{\bar{2}}\right) dz_{1}\wedge dz_{2}\wedge d\bar{z}_{2} \\
& +\left(\partial_{1}\bar{\partial}_{2}f_{1}-\bar{\partial}_{2}^{2}f_{2}+\bar{\partial}_{1}\bar{\partial}_{2}f_{\bar{1}}-\partial_{2}\bar{\partial}_{2}f_{\bar{2}}\right) dz_{1}\wedge d\bar{z}_{1}\wedge d\bar{z}_{2} \\
& +\left(\partial_{1}^{2}f_{1}-\partial_{1}\bar{\partial}_{2}f_{2}+\partial_{1}\bar{\partial}_{1}f_{\bar{1}}-\partial_{1}\partial_{2}f_{\bar{2}}\right) dz_{2}\wedge d\bar{z}_{1}\wedge d\bar{z}_{2}.
\end{align*}
To satisfy the source-free Maxwell's equations, we need $d\star F_{\omega}=0$, which in matrix form is
\[\text{diag}\{\partial_{2},\bar{\partial}_ {1},\bar{\partial}_{2},\partial_{1}\}
\begin{pmatrix}
\partial_{1} & -\bar{\partial}_{2} & \bar{\partial}_{1} & -\partial_{2}\\
\partial_{1} & -\bar{\partial}_{2} & \bar{\partial}_{1} & -\partial_{2}\\
\partial_{1} & -\bar{\partial}_{2} & \bar{\partial}_{1} & -\partial_{2}\\
\partial_{1} & -\bar{\partial}_{2} & \bar{\partial}_{1} & -\partial_{2}
\end{pmatrix}
\begin{pmatrix} f_{1} \\ f_{2} \\ f_{\bar{1}} \\ f_{\bar{2}}
\end{pmatrix}= \mathbf{0}.
\]
Clearly, this is true if and only if $\partial_{1}f_{1}-\bar{\partial}_{2}f_{2}+\bar{\partial}_{1}f_{\bar{1}}-\partial_{2}f_{\bar {2}}$ is constant and this completes the proof.\qed \\

\section{On the Lorenz gauge}

It was indicated in Section 3.2 that given a smooth $4$-vector $(\phi, A_1, A_2, A_3)$ one can associate with it the magnetic potential $1$-form $\omega=\phi dx_{0}-A_{1}dx_{1}-A_{2}dx_{2}-A_{3}dx_{3}.$ The
Lorenz gauge condition (\ref{eq:3.12}) stipulates the normalization (\ref{eq:3.15}) regarding the sum of partial derivatives, namely,
\begin{equation*}
\frac{\partial\phi}{\partial x_{0}}+\frac{\partial A_{1}}{\partial x_{1}}+\frac{\partial A_{2}}{\partial x_{2}}+\frac{\partial A_{3}}{\partial x_{3}}=0.
\end{equation*}
Theorems \ref{thm:main1} and \ref{thm:main2} indeed give a mathematical explanation as to why the Lorenz gauge matters. Here we give a unified treatment.

\begin{cor}
Let $\omega$ be a smooth $1$-form as defined in \textup{(\ref{eq:3.20})}. Then the sum of partial derivatives appearing in 
\textup{Theorem \ref{thm:main1}} is $-\frac{1}{2}d^*\omega$, and that in \textup{Theorem \ref{thm:main2}} is $\frac{1}{2}d^*\omega$.
\end{cor}
\proof
First, recall that in the Euclidean metric over $\mathbb{C}^{2}$ we have that $d^{*}\omega =-\star d \star \omega$. 
Then using the calculations in Example \ref{duals} one easily verifies that
\[
\star\omega= \frac{1}{2}\left(f_{1}dz_{1}\wedge dz_{2}\wedge d\bar{z}_{2}+f_{\bar{1}}dz_{2}\wedge d\bar{z}_{1}\wedge d\bar{z}_{2}-f_{2}dz_{1}\wedge dz_{2}\wedge d\bar{z}_{1}-f_{\bar{2}}dz_{1}\wedge d\bar{z}_{1}\wedge d\bar{z}_{2}\right).
\]
It follows that 
\[\bar{\partial}_{1}f_{1}+\bar{\partial}_{2}f_{2}+\partial_{1}f_{\bar{1}}+\partial_{2}f_{\bar {2}}=-\frac{1}{2}d^{*}\omega.\]

The sums of partial derivatives appearing in (\ref{eq:3.15}) is the real variable version of that in Theorem \ref{thm:main2}.
In the Minkowski metric over $\mathbb{C}^{2}$, we have
\begin{equation}
d^*\omega = \star d \star \omega.
\label{eq:3.38}
\end{equation}
Under the Minkowski metric, using Example \ref{mink} we have 
\begin{align}
\star\omega = \frac{1}{2} &( f_{1}dz_{2}\wedge d\bar{z}_{1}\wedge d\bar{z}_{2}+f_{2}dz_{1}\wedge dz_{2}\wedge d\bar{z}_{1}+f_{\bar{1}}dz_{1}\wedge dz_{2}\wedge d\bar{z}_{2} \label{eq:3.39} \\
 &+ f_{\bar{2}}dz_{1}\wedge d\bar{z}_{1}\wedge d\bar{z}_{2}). \nonumber
\end{align}
Now applying the exterior derivative to (\ref{eq:3.39}), after simplifying and rearranging terms, we have
\begin{align}
d\star\omega &= \frac{1}{2}\left(\partial_{1}f_{1}-\bar{\partial}_{2}f_{2}+\bar{\partial}_{1}f_{\bar{1}}-\partial_{2}f_{\bar{2}}\right)\left(dz_{1}\wedge d\bar{z}_{1}\wedge dz_{2}\wedge d\bar{z}_{2}\right) \label{eq:3.40} \\
&= -2\left(\partial_{1}f_{1}-\bar{\partial}_{2}f_{2}+\bar{\partial}_{1}f_{\bar{1}}-\partial_{2}f_{\bar{2}}\right)\left(dx_{0}\wedge dx_{1}\wedge dx_{2}\wedge dx_{3}\right). \nonumber
\end{align}
Since $\star\left(dx_{0}\wedge dx_{1}\wedge dx_{2}\wedge dx_{3}\right)=-1$ in the Minkowski metric, it follows that
\begin{equation}
d^{*}\omega= 2\left(\partial_{1}f_{1}-\bar{\partial}_{2}f_{2}+\bar{\partial}_{1}f_{\bar{1}}-\partial_{2}f_{\bar{2}}\right).
\label{eq:3.41}
\end{equation}
\qed\\

In the case we write $\omega$ in the real form $\phi dx_{0}-A_{1}dx_{1}-A_{2}dx_{2}-A_{3}dx_{3},$ then 
(\ref{eq:3.41}) implies 
\[d^*\omega=-\left(\frac{\partial\phi}{\partial x_{0}}+\frac{\partial A_{1}}{\partial x_{1}}+\frac{\partial A_{2}}{\partial x_{2}}+\frac{\partial A_{3}}{\partial x_{3}}\right).\]

If $d^{*}\omega$ is a constant, say $k$, then one can easily modify $f_i, f_{\bar{i}}, i=1, 2$ such that $d^{*}\omega$ becomes $0$. For example, in the Euclidean metric we can replace $f_1$ by $f_1+\frac{k}{2}\bar{z}_1$ and keep other functions unchanged. Similar modification can be done in the Minkowski metric. In this view, the Lorenz gauge condition is just a trivial strengthening of the condition $d^{*}\omega$ being constant. Therefore, we shall say that a smooth $1$-form $\omega$ satisfies the {\em Lorenz gauge condition} if $d^*\omega$ is constant. 

With the foregoing observation, in the Euclidean metric case one can write (\ref{eq:3.23}) as 
\begin{align*}
d\star F_{\omega} & =\left(2\nabla^{2}f_{2}+\frac{1}{2}\partial_{2}d^*\omega \right)dz_{1}\wedge dz_{2}\wedge d\bar{z}_{1}
 -\left( 2\nabla^{2}f_{1}+\frac{1}{2}\partial_{1}d^*\omega   \right)dz_{1}\wedge dz_{2}\wedge d\bar{z}_{2}\\
 & +\left( 2\nabla^{2}f_{\bar{2}}+\frac{1}{2}\bar{\partial_{2}}d^*\omega  \right)dz_{1}\wedge d\bar{z}_{1}\wedge d\bar{z}_{2} -\left(2\nabla^{2}f_{\bar{1}}+\frac{1}{2}\bar{\partial_{1}}d^*\omega   \right)dz_{2}\wedge d\bar{z}_{1}\wedge d\bar{z}_{2}. 
\end{align*}
And in the Minkowski metric case, one can write (\ref{eq:3.36}) as 
\begin{align*}
2d\star F_{\omega} & =\left(-\Box^{2}f_{2}+\partial_{2}d^*\omega \right)dz_{1}\wedge dz_{2}\wedge d\bar{z}_{1} +\left(- \Box^{2}f_{1}+\bar{\partial_{1}}d^*\omega   \right)dz_{1}\wedge dz_{2}\wedge d\bar{z}_{2}\\
 & +\left( -\Box^{2}f_{\bar{2}}+\bar{\partial_{2}}d^*\omega  \right)dz_{1}\wedge d\bar{z}_{1}\wedge d\bar{z}_{2} +\left(-\Box^{2}f_{\bar{1}}+{\partial_{1}}d^*\omega   \right)dz_{2}\wedge d\bar{z}_{1}\wedge d\bar{z}_{2}. 
\end{align*}
Moreover, if $\omega$ is a solution to the Maxwell's equations, we have $d^{*}d\omega=d^*F_\omega=0$. Hence $\Delta\omega =(dd^{*}+d^{*}d)\omega=0$ if and only if $dd^*\omega=0$, i.e., $d^*\omega$ is a constant, or in other words $\omega$ satisfies the Lorenz gauge condition. We summarize Theorem \ref{thm:main1}, Theorem \ref{thm:main2}, and the foregoing observations in the following corollary.
\begin{cor} \label{omega}
Let $\omega$ be a solution to the Maxwell's equations in vacuum under the Euclidean or Minkowski metric. Then the following are equivalent.
\begin{enumerate}[label=\textup{(\alph*)},leftmargin=*] 
\item $\omega$ satisfies the Lorenz gauge condition.
\item $\omega$ is harmonic or respectively wavelike.
\item $\omega$ is HL-harmonic.
\end{enumerate}
\end{cor}

If $\omega$ as defined in \textup{(\ref{eq:3.20})} is a solution to the Maxwell's equations in vacuum under the Euclidean or Minkowski metric, then for every smooth function $u$ the form $\omega'=\omega+du$ is also a solution because 
\[F_{\omega'}=d\omega'=d\omega+d^2u=d\omega=F_\omega.\]
This is the gauge invariance of the Maxwell's equations in differential forms. We write the above gauge transformations as 
\[f'_j=f_j+\partial_j u,\ \ \ f'_{\bar{j}}=f_{\bar{j}}+\bar{\partial}_ju,\ j=1, 2.\]
Then with respect to Theorem \ref{thm:main1} direct computations give
\begin{equation}\label{gauge1}
\bar{\partial}_1f'_1+\bar{\partial}_2f'_2+\partial_1f'_{\bar{1}}+\partial_2f'_{\bar{2}}=
\bar{\partial}_1f_1+\bar{\partial}_2f_2+\partial_1f_{\bar{1}}+\partial_2f_{\bar{2}}+\frac{1}{2}\nabla^2u,
\end{equation}
and likewise with respect to Theorem \ref{thm:main2} we have
\begin{equation}\label{gauge2}
{\partial}_1f'_1-\bar{\partial}_2f'_2+\bar{\partial}_1f'_{\bar{1}}-\partial_2f'_{\bar{2}}=
{\partial}_1f_1-\bar{\partial}_2f_2+\bar{\partial}_1f_{\bar{1}}-\partial_2f_{\bar{2}}+\frac{1}{2}\Box u.
\end{equation}
It is known that (\cite{evans2010partial, hassani2013mathematical}) for every smooth function $h$ on ${\mathbb R}^4$, the equations
\[\nabla^2u=h,\ \ \text{and}\ \  \Box u=h\]
both have solutions (non-unique). Hence there exists a smooth function $u$ such that
$\omega'=\omega+du$ satisfies the Lorenz gauge condition with respect to the Euclidean metric (or the Minkowski metric). Hence by Corollary \ref{omega} the curvature field $F_\omega=F_{\omega'}=d\omega'$ is harmonic (or respectively wavelike).
We summarize this observation as follows.

\begin{cor}
Let $F_\omega$ be a solution to the Maxwell's equations in vacuum under the Euclidean or Minkowski metric. Then
$F_\omega$ is harmonic, or respectively wavelike.
\end{cor}

In particular, this indicates that there is no non-wavelike solution to the Maxwell's equations in vacuum.

\section{Concluding Remarks}

Complex analysis is a core component in mathematics, and it has also played an increasingly important role in modern physics. It is thus meaningful to reinterpret some fundamental theories in physics from a complex perspective, for instance special relativity, Maxwell's equations, and Yang-Mills equations, whose original fomulations were in real variables. This reinterpretation will not only provide a complex formulation of the theories, but also give rise to new and natural observations from this point of view. The exploration in this direction has been made in literature, see for example \cite{hoyos2015new, kleefeld2012complex}, but it is far from being complete. This paper shall serve as a starting point for the authors to explore greater applications of complex analysis to physics theories.\\

{\bf Acknowledgments.} The authors would like to thank Marius Beceanu and Oleg Lunin for valuable comments on the initial draft of this paper. This paper is in part based on the first author's doctoral dissertation (\cite{munshi2020maxwell}) submitted to SUNY at Albany, and he is grateful to the Department of Mathematics and Statistics for providing him an opportunity to pursue his research interests.

\end{document}